\documentclass{article}

\usepackage{amsfonts}
\usepackage{amsmath}
\usepackage{amssymb}
\usepackage{amsthm}
\usepackage{amscd}
\usepackage{verbatim}
\usepackage[dvips]{graphicx}
\usepackage{graphics}
\usepackage[all,2cell]{xypic}

\theoremstyle{plain}
\newtheorem{thm}{Theorem}[section]
\newtheorem{con}[thm]{Conjecture}
\newtheorem{cor}[thm]{Corollary}
\newtheorem{lem}[thm]{Lemma}
\newtheorem{prop}[thm]{Proposition}
\newtheorem{que}[thm]{Question}

\theoremstyle{definition}
\newtheorem{cau}[thm]{Caution}
\newtheorem{conc}[thm]{Conclusion}
\newtheorem{df}[thm]{Definition}
\newtheorem{ex}[thm]{Example}
\newtheorem{exc}[thm]{Excuse}
\newtheorem{his}[thm]{Historical Note}
\newtheorem{nt}[thm]{Notations}
\newtheorem{ob}[thm]{Observations}
\newtheorem{rem}[thm]{Remark}
\newtheorem{rev}[thm]{Review}
\newtheorem{sol}[thm]{Solution}
\newtheorem{var}[thm]{Variant}

\theoremstyle{remark}

\DeclareMathOperator{\id}{id}
\DeclareMathOperator{\isoto}{\overset{\scriptstyle{\sim}}{\to}}
\DeclareMathOperator{\pt}{pt}

\DeclareMathOperator{\Ker}{Ker}

\DeclareMathOperator{\kos}{Kos}

\DeclareMathOperator{\Codim}{Codim}
\DeclareMathOperator{\Spec}{Spec}
\DeclareMathOperator{\Supp}{Supp}
\DeclareMathOperator{\Perf}{Perf}
\DeclareMathOperator{\hight}{ht}

\DeclareMathOperator{\Ar}{Ar}
\DeclareMathOperator{\Ob}{Ob}

\DeclareMathOperator{\op}{op}

\DeclareMathOperator{\Hom}{Hom}
\DeclareMathOperator{\HOM}{\mathcal{HOM}}
\DeclareMathOperator{\HHom}{\mathbf{Hom}}

\DeclareMathOperator{\Lax}{Lax}
\DeclareMathOperator{\LAX}{\mathcal{LAX}}
\DeclareMathOperator{\Aut}{Aut}

\DeclareMathOperator{\diag}{diag}
\DeclareMathOperator{\colim}{colim}
\DeclareMathOperator{\hocolim}{hocolim}

\DeclareMathOperator{\B}{B}

\DeclareMathOperator{\N}{N}
\DeclareMathOperator{\Homo}{H}

\DeclareMathOperator{\AAA}{\mathcal{A}}

\DeclareMathOperator{\EEE}{\mathcal{E}}
\DeclareMathOperator{\FFF}{\mathcal{F}}
\DeclareMathOperator{\HHH}{\mathcal{H}}

\DeclareMathOperator{\MMM}{\mathcal{M}}
\DeclareMathOperator{\PPP}{\mathcal{P}}
\DeclareMathOperator{\SSS}{\mathcal{S}}
\DeclareMathOperator{\VVV}{\mathcal{V}}
\DeclareMathOperator{\WWW}{\mathcal{W}}
\DeclareMathOperator{\XXX}{\mathcal{X}}
\DeclareMathOperator{\YYY}{\mathcal{Y}}

\DeclareMathOperator{\Set}{\mathbf{Set}}
\DeclareMathOperator{\Top}{\mathbf{Top}}
\DeclareMathOperator{\Cat}{\mathbf{Cat}}
\DeclareMathOperator{\delset}{\Delta^{op}\mathbf{Set}}
\DeclareMathOperator{\delcat}{\Delta^{op}\mathbf{Cat}}
\DeclareMathOperator{\ulG}{\underline{G}}
\DeclareMathOperator{\Fib}{\textbf{Fib}}
\DeclareMathOperator{\Kos}{\textbf{Kos}}

\DeclareMathOperator{\Sh}{Sh}
\DeclareMathOperator{\ev}{ev}

\DeclareMathOperator{\er}{er}
\DeclareMathOperator{\fib}{fib}
\DeclareMathOperator{\cof}{cof}

\DeclareMathOperator{\lax}{lax}
\DeclareMathOperator{\qis}{qis}
\DeclareMathOperator{\strict}{s}
\DeclareMathOperator{\rel}{rel}
\DeclareMathOperator{\edge}{e}
\DeclareMathOperator{\ori}{ori}

\UseAllTwocells

\title{Gersten's conjecture for commutative discrete valuation rings}
\author{Satoshi Mochizuki
\footnote{This research is supported by the 21 century COE program at Graduate School of Mathematical Sciences, the University of Tokyo.
}}
\date{}

\begin{document}

\maketitle

\vspace{-1cm}
\begin{center}
\tt{mochi@ms.u-tokyo.ac.jp}
\end{center}

\begin{abstract}
The purpose of this article is to prove that Gersten's conjecture for a 
commutative discrete valuation ring is true. 
Combining with the result of \cite{GL87}, we 
learn that Gersten's conjecture is true if the ring is a commutative regular local, smooth over a commutative discrete valuation ring.  
\end{abstract}

\tableofcontents

\section{Introduction}

In \cite{Ger73} p.28, Gersten proposed the following conjecture

\begin{con} [Gersten's conjecture] \label{GC} $ $
\\
For any commutative regular local ring $A$ and natural numbers $n$, $p$, the canonical inclusion 
$\MMM^{p+1}(A) \hookrightarrow \MMM^p(A)$ induces the zero 
map on $K$-groups
$$K_n(\MMM^{p+1}(A)) \to K_n(\MMM^p(A)) \\\ ,$$
where $\MMM^i(A)$ is the category of finitely generated $A$-modules $M$
with\\
 $\Codim_{\Spec A} \Supp M \geqq i$. 
\end{con}

\begin{his} \label{HN}$ $
\\
The conjecture has been proved in the following cases:\\
(1) {\textbf{The case of general dimension}}\\
(a) If $A$ is of equi-characteristic, 
then Gersten's conjecture for $A$ is true.\\
We refer to \cite{Qui73} for special cases, and the general cases \cite{Pan03} can be deduced from limit argument and Popescu's general N\'eron desingularization 
\cite{Pop86}. (For the case of commutative discrete valuation rings, it was first proved by Sherman \cite{She78})\\
(b) If $A$ is smooth over some commutative discrete valuation ring $S$ and satisfies some condition, then Gersten's conjecture for $A$ is true. \cite{Blo86}\\
(c) If $A$ is smooth over some commutative discrete valuation ring $S$ and if we accept Gersten's conjecture for $S$, then Gersten's conjecture for $A$ is true. \cite{GL87}\\
(2) {\textbf{The case of that $A$ is a commutative discrete valuation ring}}\\
(d) For the cases of $n=0,\ 1$ are classical and $n=2$ was proved by Dennis and Stein, announced in \cite{DS72} and proved in \cite{DS75}.\\ 
(e) If its residue field is algebraic over a finite field, then Gersten's conjecture for $A$ is true.\\
We refer to \cite{Ger73} for the cases of that its residue field is a finite field, the proof of general cases \cite{She82} is improved from that of special cases by using Swan's result \cite{Swa63}, universal property of algebraic $K$-theory \cite{Hil81} and limit argument.\\
(3) {\textbf{Other works}}\\
(f) For torsion coefficient $K$-theory, Gersten's conjecture for a commutative discrete valuation ring is true. \cite{Gil86}, \cite{GL00}\\
(g) Some conjectures imply that Gersten's conjecture for a commutative discrete valuation ring is true. \cite{She89} 
\end{his}

\noindent
The purpose of this article is to prove the following theorem.

\begin{thm} \label{main} $ $
\\
For any commutative discrete valuation rings, Gersten's conjecture is true. 
\end{thm}

\noindent
Combining \textbf{Historical Note \ref{HN}} (c) with \textbf{Theorem \ref{main}}, we have the following result.

\begin{cor} \label{cormainthm} $ $
\\
For any commutative regular local ring which is smooth over a commutative
discrete valuation ring, Gersten's conjecture is true.
\end{cor}

\begin{ex} [Due to Kazuya Kato] \label{noncom} $ $
\\
In \textbf{Theorem \ref{main}}, the assumption of commutativity is essential. (See the proof of \textbf{Lemma \ref{key lemma 2}}.)
Let $D$ be a skew field finite over
$\mathbf{Q}_p$, $A$ its integer ring and $a$ its prime element. As the inner
automorphism of $a$ induces non-trivial automorphism on its residue field,
we have $x \in A^{\times}$ with $y=axa^{-1}x^{-1}$ is non vanishing in its residue field, a fortiori in $K_1(A)={(A^{\times})}^{ab}$. On the other hand $y$ is a
commutator in $D^\times$. Hence we learn $K_1(A) \to K_1(D)={(D^{\times})}^{ab}$ is not injective.
\end{ex}

\begin{nt} $ $
\\
In this paper, we will use the following notations.\\
(1) $R$ denote a commutative discrete valuation ring, and $\pi$ its uniformizer element, $k=R/\pi$.\\ 
(2) A topological space means a compactly generated topological space and their fiber product is considered in the category of compactly generated topological spaces. A pointed space means a pointed topological space with non-degenerate base point.\\
(3) When we quote the paper \cite{Qui73}, we follow the page numbering by the numbering at the bottom. (In \cite{Qui73}, there are numbers both at the top and the bottom of each page, and the numberings do not agree.)
\end{nt}

\noindent
\textbf{Acknowledgement} Throughout the last six years, the author is helped by many people. He is greatful for them. He learned fundamental knowledge about several filtrations on algebraic $K$-theory from Kei Hagihara and about several constructions of the Adams operations from Bernhard K\"ock. Arguing with Masana Harada, he polished his techniques on algebraic $K$-theory. Kazuya Kato taught him many not only painful but also deeply impressive suggestions. Discussing the structure of Grothendieck groups of group representations over fields with Shushi Harashita and Shu Kato, the structure of strictly perfect complexes over regular local rings with Ken Naganuma and Miho Aoki, and the structure of homotopy colimits in general $2$-categories with Go Yamashita, weight argument of the Adams operations with Pierre Deligne and Christophe Soul\'e, he had been changing his delusional idea into a real mathematical proof. Takeshi Saito, Seidai Yasuda, Kentaro Nakamura, Yoshihiro Sawano and Ofer Gabber carefully read a draft version of this paper. Charles A. Weibel remarked historical comments. Kanetomo Sato, Takao Yamazaki, Fumitoshi Sato and Masaki Hanamura invited him to their universities. Kensuke Itakura, Toshiro Hiranouchi and Yuki Kato took part in his lectures on Gersten's conjecture and gave many worthy questions and led him to simplifying the proof of the main theorem in this paper. Shuji Saito have been encouraging him to complete this work. 

\section{The strategy of proving the main theorem}

The proof of Gersten's conjecture in this paper requires a long, sustained argument. In this section, we will explain how to prove the main theorem and philosophy behind the proof. The author's proof of Gersten's conjecture is inspired from the works related about Serre's conjecture for intersection multiplicities \cite{Ser65} by H. Gillet, M. Levine, and C. Soul\'e 
\cite{GS87}, \cite{GS99} and \cite{Lev85}. Their studying comparison between coniveau filtrations and $\gamma$-filtrations led him to the idea that \lq\lq weight argument of the Adams operations implies Gersten's conjecture.\rq\rq \ Now, we explain more precisely. First we will review the following \textbf{Proposition \ref{Levinelemma}}.\\

\begin{nt} $ $
\\
As in the Conjecture \ref{GC}, for noetherian finite Krull dimensional commutative ring with unit $A$, let $\MMM^p(A)$ be the category of finitely generated $A$-modules $M$ with  $\Codim_{\Spec A} \Supp M \geqq p$.\\
\end{nt}

\begin{prop} \label{Levinelemma} $ $
\\
{\rm (c.f. \cite{Lev85} P.452, Proposition 1.1 and \cite{GS87} p.263 Proposition 4.12)} Let $A$ be a regular domain.\\
{\rm (1)} The following statements are equivalent.\\ 
{\rm (i)} The maps $K_0(\MMM^p(A)) \to K_0(\MMM^{p-1}(A))$ are zero for $p=1,\cdots,\dim A$.\\
{\rm (ii)} $K_0(\MMM^p(A))$ is generated by cyclic modules $A/(f_1,\cdots,f_p)$ where $f_1,\cdots,f_p$ forms a regular sequence for $p=1,\cdots,\dim A$.\\
{\rm (2)} Moreover we consider the following statements\\
{\rm (i)$_\mathbb{Q}$} The maps $K_0(\MMM^p(A))\otimes \mathbb{Q} \to K_0(\MMM^{p-1}(A)) \otimes \mathbb{Q}$ are zero for $p=1,\cdots,\dim A$.\\
{\rm (iii)} For non-negative integer $k$, the Adams operation $\psi_k$ acts on $K_0(\MMM^p(A))$ by multiplication of $k^{p}$.\\
Then {\rm (ii)} implies {\rm (iii)} and {\rm (iii)} implies {\rm (i)$_\mathbb{Q}$}.
\end{prop}

\begin{proof}[Sketch of proof]
(1) The implication (i) $\Rightarrow$ (ii) is proved in \cite{Lev85} Ibid by using the following surjections
$$K_1(\MMM^{p-1}(A)/\MMM^p(A))=\underset{\substack{x \in \Spec A \\ \hight x=p-1}}{\bigoplus}k(x)^{\times} \to K_0(\MMM^p(A))$$
for $p=1,\cdots,\dim A$ and induction argument. The implication (ii) $\Rightarrow$ (i) can be proved by the following exact sequence and additivity of $K$-theory.
\begin{equation} \label{crutial exact seq 1}
0 \to A/(f_1,\cdots,f_{p-1}) \overset{f_p}{\to} A/(f_1,\cdots,f_{p-1}) \to   
A/(f_1,\cdots,f_p) \to 0
\end{equation}
where $f_1,\cdots,f_p$ are regular sequence of $A$. (Compare the proof of \textbf{Lemma \ref{key lemma 2}}.)\\
(2) The implication (iii) $\Rightarrow$ (i)$_\mathbb{Q}$ can be proved by a weight argument of Adams operations.\\
The implication (ii) $\Rightarrow$ (iii) will be proved as follows. First, for any closed set $Y \subset \Spec A$, we put $\Perf^{Y}(\Spec A)$ the category of strictly perfect complexes which are acyclic on $\Spec A - Y$, that is, whose any object is a bounded complex $X$ of finitely generated projective $A$-modules such that $X|_{\Spec A -Y}$ is acyclic. We also put $\Perf^p(\Spec A):=\underset{\substack{Y \subset X \\ \Codim Y \geqq p}}{\cup} \Perf^Y(\Spec A)$.\\
\\
\noindent
\textbf{Claim}\\
For any non-negative integer $n$, we have the following identity.
$$K_n(\MMM^p(A)) \isoto K_n(\Perf^p(\Spec A);\qis)$$  
where $\qis$ means quasi-isomorphisms, that is, $K_n(\Perf^p(\Spec A);\qis)$ is a Waldhausen quasi-isomorphism \lq Klassen\rq \ group.\\
\begin{proof}[Proof of \textbf{Claim}] We have the following identities 
$$K_n(\MMM^p(A)) \underset{\text{I}}{\isoto} \underset{\substack{Y \subset \Spec A \\ \Codim Y \geqq p}}{\colim} K'_n(Y) \underset{\text{II}}{\isoto} \underset{\substack{Y \subset \Spec A \\ \Codim Y \geqq p}}{\colim} K_n(\Spec A \ \text{on}\ Y) \underset{\text{I}}{\isoto} K_n(\Perf^p(\Spec A);\qis)$$   
where isomorphisms I are proved by continuity \cite{Qui73}, \cite{TT90} and isomorphism II is proved by the Poincar\'e duality and comparing the following fibration sequences \cite{Qui73} and \cite{TT90}
$$K'(Y) \to K'(\Spec A) \to K'(\Spec A -Y) ,$$
$$K(\Spec A \ \text{on} \ Y) \to K(\Spec A) \to K(\Spec A -Y) $$
for any closed set $Y \subset \Spec A$.
\end{proof}
For $n=0$, through the isomorphism in \textbf{Claim}, a class of  $A$-module $[A/(f_1,\cdots,f_p)]$ correspond to a class of the Koszul complex $[\kos(f_1,\cdots,f_p)]$ where $f_1,\cdots,f_p$ forms a regular sequence. So we learn (ii) and (iii) are equivalent to the following statement (ii)' and (iii)' respectively.\\
(ii)' $K_0(\Perf^p(A);\qis)$ is generated by Koszul complexes $\kos(f_1,\cdots,f_p)$ where $f_1,\cdots,f_p$ forms a regular sequence for $p=1,\cdots,\dim A$.\\
(iii)' For non-negative integer $k$, the Adams operation $\psi_k$ acts on $K_0(\Perf^p(\Spec A);\qis)$ by multiplication of $k^{p}$.\\
Similar argument in \cite{GS87}, we can define the Adams operations $\psi_k$ on $K_0(\Perf^p(\Spec A);\qis)$ by using the Dold-Puppe correspondence and the Serre's result about Grothendieck groups of reductive group schemes \cite{Ser68}, and we have the identity 
$$\psi_k([\kos(f_1,\cdots,f_p)])=k^p[\kos(f_1,\cdots,f_p)]$$
by \cite{GS87} Ibid. Hence assertion (ii)' implies assertion (iii)'.
\end{proof}

\noindent
To prove Gersten's conjecture, we shall extend the argument above to a higher analogue of that.
The author's starting point to prove the Gersten's conjecture for commutative discrete valuation rings is the following computational result \textbf{Theorem \ref{comparison thm}}.

\begin{nt} $ $
\\
Let $\Kos^1(R)$ be the full additive subcategory of complexes $X$
in $\MMM(R)$ such that $X_i=0$ for $i \ne 0,1$ and $X_i$ are free $R$ modules for $i=0, 1$ and its boundary morphism $d^X:X_1 \to X_0$ is injective and homology groups are in $\MMM^1(R)$.\\
\end{nt}

\begin{rem} $ $
\\
We say that a sequence $X \to Y \to Z$ in $\Kos^1(R)$
is an admissible short exact sequence if it is a short exact sequence
in $C_b(\MMM(R))$. Then $\Kos^1(R)$ is an exact category and
$$\Homo_0:\Kos^1(R) \to \MMM^1(R)$$
is an exact functor.
\end{rem}

\begin{thm}[\cite{Moc07}] \label{comparison thm} $ $
\\
In the notation above, for any non-negative integer $n$, $\Homo_0$ induces an isomorphism
$$K_n(\Homo_0):K_n(\Kos^1(R);\qis) \isoto K_n(\MMM^1(R))$$
where $\qis$ means quasi-isomorphisms, that is, $K_n(\Kos^1(R);\qis)$ is a quasi-isomorphisms Waldhausen \lq Klassen\rq\ group.
\end{thm}

\noindent
This theorem means that higher $K$-groups of $\MMM^1(R)$ are also generated by Koszul complexes. We shall ask the following question.\\

\begin{que} $ $
\\
For any regular local ring $A$, For any non-negative integer $n$ and $p$, is $K_n(\MMM^p(A))$ generated by weight $p$ pure objects? Here the phrase \lq weight $p$ pure objects\rq \ means Koszul complexes $\kos(f_1,\cdots,f_p)$ (resp. cyclic modules $A/(f_1,\cdots,f_p)$) where $f_1,\cdots,f_p$ forms a regular sequence.
\end{que}

\noindent
\begin{ob} \label{observe} $ $
\\
(1) Now remeber that there are several methods of constructing the Adams operations on higher algebraic $K$-theory for schemes.\\
(A) Categorical discription of the Adamas operations \cite{Gra92}, \cite{Sch87}. (See also \cite{Tho77} p.111.)\\
(B) First using group representations theory, we define them for affine cases \cite{Kra80} \cite{Hil81}, and by patching argument or Jouanolou's device, generalize to it. \cite{Sou85} \cite{Lec98} \cite{Lev97}.\\
(C) Stable homotopy theoretic approach. \cite{Rio06}.\\
(2) Now consider what we need to do weight argument integrally.\\ 
(a) We don't need to precisely compute the action of $\psi_k$ on $K_n(\MMM^p(A))$, we only verify that the weight is changing thorugh the map $K_n(\MMM^{p+1}(A)) \to K_n(\MMM^p(A))$.\\
(b) This means that we shall only do the argument in the proof of \textbf{Proposition \ref{Levinelemma}} (ii) implies (i). (Compare the Quillen's proof of Gersten's conjecture \cite{Qui73} p.133-134.) So we shall only use the higher analogue of short exact sequence (\ref{crutial exact seq 1}). (Compre the short exact sequence (\ref{crutial exact seq 2}) in the proof of \textbf{Lemma \ref{key lemma 2}}.)\\  
\end{ob}

\noindent
In this paper we will choose the method (B) and do argument (b) in \textbf{Observation \ref{observe}}. Next we will explain the difficulity of doing argument (b).\\

\begin{ob} \label{observe 2} $ $
\\
It is well-known that the category of small categories, simplicial sets and (compactly generated) topological spaces are related the nerve functor \cite{Gro59} and the geometric realization functor \cite{Mil57}
$$\Cat \overset{\N}{\to} \delset \overset{|?|}{\to} \Top$$
and they induce categories equivalence between their homotopy categories respectively
\begin{equation} \label{dogma}
\SSS^{-1}\Cat \isoto \SSS^{-1}\delset \isoto \SSS^{-1} \Top
\end{equation}
where the morphism sets $\SSS$ are certain classes of weak equivalences respectively. (See the former equivalence \cite{Ilu72}, \cite{Tho77} and the latter one \cite{Qui67}.) In theoretically, we can approximate to topological spaces by combinatorial data (that is,  simplicial sets) and by special type of simplicial sets (for example, categories) up to homotopies. This thought is a fundamental dogma in algebraic $K$-theory. But in practice, there are a few morphisms between categories as in the following \textbf{Example \ref{a few mor ex}}. So many authors attempt to replace $\Cat$ with a category of more higher categorical objects. (c.f. \cite{Tho77}, \cite{Tho79}, \cite{HS85}, \cite{Sch87} and \cite{TV04}.) The author will also introduce a class of suitable $2$-categories for our argument in \S 2.\\
\end{ob}

\begin{ex} \label{a few mor ex} $ $
\\
Let $A$ be a regular noetherian ring, $\MMM(A)$ the category of finitely generated $A$-modules and $C_b(\PPP(A))$ the category of bounded complexes of finitely generated projective $A$-modules. Then we have the following assignment 
$$\Ob\MMM(A) \ni M \mapsto (\text{a projective resolution of $M$}) \in \Ob C_b(\PPP(A))$$
which induces a canonical isomorphism between two Grothendieck groups
$$K_0(\MMM(A)) \isoto K_0(C_b(\PPP(A));\qis)$$ 
where $\qis$ means a quasi isomorphisms \lq\lq Klassen\rq\rq \ group. We also have isomorphisms between their higher $K$-groups, but to prove that result is more complicated. The main reason is that the assignment above is not a functor. (Only it is a {\it{lax functor}}. See \textbf{Definition \ref{Lax functors, deformations def}} (1) and \textbf{Example \ref{lax functor ex}}.) So in the usual theory, we can not construct maps between higher $K$-groups from this assignment directly. In the usual theory, we used to prove this result by inserting the following inclusion functors 
$$\MMM(A) \hookleftarrow \PPP(A) \hookrightarrow C_b(\PPP(A))$$ 
and verifying that they induce isomorphisms on their higher algebraic $K$-groups
$$K_n(\MMM(A))\overset{\sim}{\leftarrow} K_n(\PPP(A)) \isoto K_n(C_b(\PPP(A));\qis)$$
for any non-negative integer $n$. (The former isomorphism can be proved by resolution theorem \cite{Qui73} and the latter one can be done by Gillet-Waldhausen theorem \cite{Gil81}, \cite{Wal85}, \cite{TT90}.)\\
This is just another typical example of lacking techniques on inducing morphisms between higher algebraic $K$-theory from exotic functors. (But see \cite{Nee00}.)\\
If we discuss about using sequence (\ref{crutial exact seq 1}), we shall consider the following assignment
$$\MMM^p(A) \to C_b(\MMM^{p-1}(A))$$ 
$$R/(f_1,\cdots,f_p) \mapsto 
\begin{bmatrix}
R/(f_1,\cdots,f_{p-1})\\
\downarrow {f_p}\\
R/(f_1,\cdots,f_{p-1})
\end{bmatrix}$$
In general, this is not functorial in usual sense. (Compare the proof of \textbf{Lemma \ref{key lemma 2}}, especially the short exact sequence (\ref{crutial exact seq 2}))\\
\end{ex}

\begin{sol}[Rectification techniques] \label{rectification}$ $
\\
Here is one of the solution of the problem above, called {\it{rectification}} techniques. Now we will explain it.\\
(1) Let $f:\XXX \to \YYY$ be a functor between small categories. we can consider $\XXX$ as a family over $\YYY$. So if $\YYY$ is enough closed under taking colimits, we can reconstruct $f$ by gathering each fiber. That is, we have a formula
$$f(x)=\colim f|x(=\underset{y \to x \in \Ob\XXX/x}{\colim} f(y)).$$
(2) Let $f:\XXX \to \YYY$ be a {\it{lax functor}} in some sense and assume $\YYY$ is enough closed under taking {\it{homotopy colimits}} in some sense. Then $f$ can be rectified by the functor $\tilde{f}$ defined by the following formula
$$\tilde{f}(x)=\hocolim f|x(=\underset{y \to x \in \Ob\XXX/x}{\hocolim} f(y)).$$
See for example \cite{GR71} Expos\'e VI \S 8, \cite{Tho77} p.37 II \S 4, \cite{Tho82} p.1646-1647.\\
(3) So it is important to examine the behaviour of homotopy colimits through the equivalences in (\ref{dogma}). (See through the intentions behind the proof of Theorem B in \cite{Qui73} p.97, see also \cite{Tho82} p.1625-1626.)
\end{sol}

\noindent 
In this paper, the author will use up dated version of rectification technique by using the result of Bullesjos, Cegarra, Moerdijk and Svensson \cite{BC03} \cite{MS93} in \S 3. (See also \textbf{Remark \ref{FPVN}} (3) and compare the argument in \textbf{Solution \ref{rectification}} (3).) Next we analyze the method (B) in \textbf{Observation \ref{observe}}.

\begin{ob} \label{observe 3} $ $
\\
Quillen had constructed the universal map using study of characteristic classes of representations. (See \cite{Qui76}, \cite{Hil81}, see also \cite{Gro68}.) In \cite{She82}, Sherman read more formal argument into Quillen's construction. (Compare \textbf{Example \ref{Sherman argument}}.) Namely the construction of the universal map is only depend on the behaviour of internal hom objects through the equivalences in (\ref{dogma}). It is an important fact that to state the universal property of algebraic $K$-theory, we shall only consider skeletons of connected groupoids (that is, groups) for the domain variable of internal hom functors. This is the reason why the universal property of algebraic $K$-theory can be related with group representations.
\end{ob}

\noindent
To combine the rectification techniques with the universal property of algebraic $K$-theory, we will introduce a notion of {\it{lax group representations}}. In \S 4-1, we will study the behaviour (lax) internal hom objects and use this results, in \S 4-2, define the lax version of the universal map, in \S 4-3, compare this universal map and the original one in the certain cases.\\

\begin{cau} $ $
\\
To use the universal property of algebraic $K$-theory associated with semi-simple exact categories, we replace $\MMM^1(R)$ with $\MMM(k)$. By d\'evissage theorem, we have group isomorphisms
$$K_n(\MMM(k))\isoto K_n(\MMM^1(R))$$ 
for any non negative integer $n$. This is only group isomorphism and is not considered as special $\lambda$-ring morphism. So we shall not use the usual Adams operations on $K_n(k)$.
\end{cau}

\section{Categories with equivalence relations}

In this section, we will introduce categories with equivalence relations - a special class of $2$-categories - and lax functors. Using these notions, we will formulate lax group representations. We will also explain how to associate categories with equivalence relations with topological spaces.\\

\begin{df}[Category with equivalence relations] \label{cat with er def} $ $
\\
(1) A {\it{category with equivalence relations}} is a category $\XXX$ which has compatible equivalence relations in each hom set. That is, for each objects $X$, $Y$ in
$\XXX$, $\Hom_{\XXX}(X,Y)$ has equivalence relation $\sim$ and satisfies the following axiom.\\
Let $a:X \to Y$, $b:Y \to Z$ and $x,\ y:Y \to Z$ be morphisms in $\XXX$ and assume $x \sim y$. Then $x \circ a \sim y \circ a$ and $b \circ x \sim b \circ y$.\\(2) A morphism between categories equivalence relations $\XXX$, $\YYY$ is a functor $f:\XXX \to \YYY$ which preserves equivalence relations. That is, if $a,\ b:X\to Y$ are morphisms in $\XXX$ and assume $a \sim b$, then we have $f(a) \sim f(b)$.\\
(3) We will denote the category of small categories with equivalence relations by $\Cat_{\er}$.  
\end{df}

\begin{ex} \label{ex of cat with er} $ $
\\
(1) Any category $\XXX$ can be considered as a category with equivalence relations where equivalence relations between morphisms $f,\ g:X \to Y$ in $\XXX$ is defined by $f \sim g$ if and only if $f=g$. we will call it {\it{trivial equivalence relations}} in $\XXX$\\
(2) Let $\XXX$ be category with equivalence relations. Then its subcategory $\YYY$ is considered as a category with equivalence relations in the natural way.\\
(3) Let $\XXX$ be category with equivalence relations and $f:\YYY \to \XXX$ a category equivalence. Then $\YYY$ can be considered as a category with equivalence relations such that $f:\XXX \to \YYY$ is those morphisms in the obvious way.\\ 
(4) Let $\XXX$ be an additive category and $C_b(\XXX)$ be the category of bounded complexes in $\XXX$. Then $C_b(\XXX)$ can be considered as a category with equivalence relations. Here $f \sim g$ if and only if $f$ and $g$ are chain homotopic.\\
(5) We define the category $\AAA$ that is full subcategory of $\MMM(R)$. This plays crucial role in the proof of the main theorem. We put
$$\Ob\AAA:=\{M\in \Ob\MMM(R);\pi T(M)=0\}$$
where $T(M):=\{x \in M; rx=0\ {\text{for some}}\ r \in R\}$. 
So every object in $\AAA$ is isomorphic to the following type module $R^{\oplus n} \oplus (R/\pi)^{\oplus m}$. (If readers worry about set theoretic problems, you shall better replace the definition of $\AAA$ by its skeletons. c.f. (3).) For any morphisms
$$a=
\begin{pmatrix}
a_{f \to f} & 0\\
a_{f \to t} & a_{t \to t}
\end{pmatrix},\ 
b=
\begin{pmatrix}
b_{f \to f} & 0\\
b_{f \to t} & b_{t \to t}
\end{pmatrix}
:R^{\oplus n} \oplus (R/\pi)^{\oplus m} \to R^{\oplus u} \oplus (R/\pi)^{\oplus v}$$
where the indexes $f$ and $t$ means \lq free\rq\ and \lq torsion\rq\ respectively, $a \sim b$ if and only if $a_{f \to f} = b_{f \to f} \mod \pi$, $a_{f \to t}=b_{f \to t}$ and $a_{t \to t}=b_{t \to t}$. This definition does not depend on the choice of isomorphisms and $\AAA$ is a category with equivalence relations.

\begin{proof}
In the notation above, we assume $a \sim b$. For any $c:R^{\oplus n'} \oplus (R/\pi)^{\oplus m'} \to R^{\oplus n} \oplus (R/\pi)^{\oplus m}$ and $d:R^{\oplus u} \oplus (R/\pi)^{\oplus v} \to R^{\oplus u'} \oplus (R/\pi)^{\oplus v'}$, we want to prove $a c \sim b c$ and $d a \sim d b$. (This fact implies independence of the choice of isomorphisms.) We will use similar index notations above for $c$ and $d$. Obviously $a_{f \to f} c_{f \to f}=b_{f \to f} c_{f \to f} \mod \pi$, $d_{f \to f} a_{f \to f}=d_{f \to f} b_{f \to f} \mod \pi$, $a_{t \to t} c_{t \to t}=b_{t \to t} c_{t \to t}$, $d_{t \to t} a_{t \to t}=d_{t \to t} b_{t \to t}$ and $a_{f \to t} c_{f \to f} + a_{t \to t} c_{t \to t}=b_{f \to t} c_{f \to f} + b_{t \to t} c_{t \to t}$. Since $a_{f \to f}-b_{f \to f}$ can be divided by $\pi$, we have $d_{f \to t}(a_{f \to f}-b_{f \to f})=0$. So we learn $d_{f \to t} a_{f \to f} + d_{t \to t} a_{f \to t}=d_{f \to t} b_{f \to f} + d_{t \to t} b_{f \to t}$.
\end{proof}
\end{ex}

\begin{df}[Lax functors, Deformations] \label{Lax functors, deformations def} $ $
\\
Let $I$ be a small category and  $\XXX$ a small category with equivalence relations.\\
(1) A {\it{lax functor}} $x:I \to \XXX$ consists of two maps which assign to each object $i \in I$ an object $x_i \in \XXX$, to each morphism $a:i \to j$ in $I$ an morphism $x_a:x_i \to x_j$ in $\XXX$, respectively. These data are required to satisfy the following conditions:\\
\textbf{LF1} for any object $i\in I$, $x_{\id_i}=\id_{x_i}$,\\
\textbf{LF2} for any pair of composable morphisms $i \overset{a}{\to} j \overset{b}{\to} k$ in $I$, $x_{ba} \sim x_bx_a$.\\
(2) A {\it{(lax) deformation}} $f:x \Rightarrow y$, between lax functors $x,\ y:I \to \XXX$ consists of a map which assign to each object $i \in I$ a morphism $f_i:x_i \to y_i$ in $\XXX$. These data satisfy the condition that for each morphism $a:i \to j$ in $I$, $f_j x_a \sim y_a f_i$.\\
(3) In the notation above, $f$ is called {\it{strict}} if each morphism $a:i \to j$ in $I$, $f_j x_a=y_a f_i$.\\
(4) Let $\Theta$ be a subset of $\Ob I$ and $x,y:I \to \XXX$ a lax functors which coincide over $\Theta$. A deformation $f:x \Rightarrow y$ is qualified as {\it{relative}} to $\Theta$ whenever $f_i=\id_{x_i}$ for all $i \in \Theta$.\\
(5) For any deformations $x \overset{f}{\Rightarrow} y \overset{g}{\Rightarrow} z$, between lax functors $x,\ y ,\ z:I \to \XXX$, we can define their composition $x \overset{gf}{\Rightarrow} z$ by $(gf)_i=g_if_i$ for any $i$ in $I$. Let  $\LAX(I,\XXX)$ (resp. $\LAX(I,\XXX)_{\strict}$) be the category whose objects are lax functors from $I$ to $\XXX$ and morphisms are deformations (resp. strict deformations). $\HOM(I,\XXX)$ (resp. $\HOM(I,\XXX)_{\strict}$) is a full subcategory of $\LAX(I,\XXX)$ (resp. $\LAX(I,\XXX)_{\strict}$) whose objects are just functors. $\LAX(I,\XXX)_{\rel\Theta}$ (resp. $\HOM(I,\XXX)_{\rel\Theta}$) is a subcategory of $\LAX(I,\XXX)$ (resp. $\HOM(I,\XXX)$) whose morphisms are deformations relative to $\Theta$. And we will put $\Lax(I,\XXX):=\Ob \LAX(I,\XXX)$.
\end{df}

\begin{ex} \label{lax functor ex} $ $
\\
In the notation \textbf{Example \ref{a few mor ex}}, we can consider $C_b(\PPP(A))$ as a category with equivalence relations by \textbf{Example \ref{ex of cat with er}} (4). Then the assignment
$$\MMM(A) \to C_b(\PPP(A))$$
$$ M \mapsto \text{(a projective resolution of $M$)}$$ 
$$ f \mapsto \text{(a lifting of $f$)}$$
defines a lax functor.
\end{ex}

\begin{rem} \label{LAX hom is category with equivalence relations} $ $
\\
(1) Let $I$ and $J$ be small categories. We consider $J$ as a category with trivial equivalence relations. Then a lax functor $I \to J$ is just a usual functor. So we have identities 
$$\LAX(I,J)=\LAX(I,J)_{\strict}=\HOM(I,J)=\HOM(I,J)_{\strict}.$$
(2) Let $I$ be a small category and $\XXX$ a small category with equivalence relations. Then $\LAX(I,\XXX)$ can be considered as a category with equivalence relations where equivalence relations between deformations defined by $f \sim g$ if and only if $f_i \sim g_i$ for any objects in $i \in I$.\\
(3) Let $f:I \to J$ be a functor between small categories, $g:\XXX \to \YYY$ a morphism between small categories with equivalence relations, $x, y:J \to \XXX$ lax functors and $a:x \Rightarrow y$ a deformation. Then we can define their compositions. That is, a lax functor $gxf:I \to \YYY$ and a deformation
 $gaf:gxf \Rightarrow gyf$ in the obvious ways.\\
(4) So we have functors 
$$\LAX(?,?)_{\square},\ \HOM(?,?)_{\square}:\Cat^{\op} \times \Cat_{\er} \to \Cat_{\er}$$
where $\square =\emptyset$, $\strict$, or $\rel\Ob(?)$. The last one means the functors
$$(I,\XXX) \mapsto \LAX(I,\XXX)_{\rel \Ob I} \ \text{(resp. $
\HOM(I,\XXX)_{\rel \Ob I})$}.$$
\end{rem}

\begin{ex}[Lax group representations] \label{Lax group representations def} $ $
\\
Let $G$ be a group and $\XXX$ be a small category with equivalence relations. As usual we will sometimes consider $G$ as the category whose object set is only one element. We will denote the \lq\lq category\rq\rq \ $G$ by $\ulG$. An object in $\LAX(\ulG,\XXX)$ is a pair $(X,\rho_X)$ consisting of an object $X$ in $\XXX$ and a map $\rho_X:G \to \Aut X$ which satisfies $\rho_X(gh) \sim \rho_X(g)\rho_X(h)$ for any $g,\ h \in G$. It is called a {\it{lax group representation}} in $\XXX$. This concept plays crucial role in the proof of main theorem.   
\end{ex}

\begin{df}[Various nerves] \label{Nerves def} $ $
\\
Let $\XXX$ be a small category with equivalence relations.\\
(1) (\cite{Dus02}, \cite{Str87}) A {\it{geometric nerve}} $\Delta \XXX$ is a simplicial set defined by
$$\Delta\XXX:=\Lax(?,\XXX):\Delta^{\op} \to \Set.$$
(2) A {\it{simplicial geometric nerve}} $\underline{\Delta}\XXX$ is a simplicial category defined by
$$\underline{\Delta}\XXX:=\LAX(?,\XXX)_{\rel\Ob(?)}:\Delta^{\op} \to \Cat.$$
(3) A {\it{2-nerve}} is a simplicial category defined by
$$\N^{\er}\XXX:=\HOM(?,\XXX)_{\rel\Ob(?)}:\Delta^{\op} \to \Cat.$$
(4) Every assignments are functorial in $\XXX$ (cf. \textbf{Remark \ref{LAX hom is category with equivalence relations}} (4)), so we get the functors
$$\Delta:\Cat_{\er} \to \delset\ \text{and} \ \underline{\Delta},\ \N^{\er}:\Cat_{\er} \to \delcat.$$ 
\end{df}

\begin{rem}[Fundamental properties of various nerves] \label{FPVN} $ $
\\
Let $\XXX$ be a small category with equivalence relations.\\
(1) If the equivalence relations in $\XXX$ are trivial, then $\Delta\XXX$ is just a usual nerve $\N\XXX$ by \textbf{Remark \ref{LAX hom is category with equivalence relations}} (1).\\ 
(2) We will define a canonical natural transformation
$$\Lax(?,?) \Rightarrow \Hom_{\delset}(\N(?),\Delta(?))$$
between functors $\Cat^{\op} \times \Cat_{\er} \to \Set$, as follows.\\
For each lax functor $x:I \to \XXX$ from a small category $I$ to a small category with equivalence relations $\XXX$, we will assign a simplicial map $\hat{x}:\N I \to \Delta \XXX$ in the following way. For each non-negative integer $n$,
$$\N_n I=\Hom_{\Cat}([n],I) \ni \phi \mapsto x \circ \phi \in \Lax([n],\XXX)=\Delta_n\XXX.$$
This define the assignment $x \mapsto \hat{x}$ and it is just the desired natural transformation.\\ 
(3) The following canonical inclusions of simplicial categories induce homotopy equivalence for their geometric realizations, respectively
$$\N^{\er}\XXX \hookrightarrow \underline{\Delta}\XXX \hookleftarrow \Delta\XXX$$
where $\Delta\XXX$ is considered as a simplicial discrete category. (See \cite{MS93}, \cite{BC03}.)
\end{rem}

\section{Lax algebraic $K$-theory}

In this section, we will define ad hoc lax algebraic $K$-theory. As explained in \S 1, to introduce this concept, especially \textbf{Theorem \ref{rectification variant}}, is one of the variant of rectification techniques on algebraic $K$-theory.

\begin{df}[Categories with fibrations and equivalence relations] \label{cat withfib and er def} $ $
\\ 
(1) A {\it{category with fibrations}} is the dual of a usual category with cofibrations. (c.f.\cite{Wal85}, \cite{TT90}, \cite{TV04}). So a category with fibrations will be a pair $(\XXX,\fib(\XXX))$ consisting of a category $\XXX$ with subcategory $\fib(\XXX) \hookrightarrow \XXX$ whose morphisms will be called {\it{fibrations}} satisfying the following axioms:\\
\textbf{Fib 0} $\XXX$ has a specific zero object $\ast$.\\
\textbf{Fib 1} The subcategory $\fib(\XXX)$ contains all isomorphisms in $\XXX$.\\
\textbf{Fib 2} For any object $x \in \XXX$, the morphism $x \to \ast$ is in $\fib(\XXX)$.\\
\textbf{Fib 3} If $x \to y$ is fibration, then, for any morphism $z \to y$ in $\XXX$, the fiber product $x \times_y z$ exists in $\XXX$ and the canonical morphism $x \times_y z \to z$ is again a fibration.\\  
We always omit $\fib(\XXX)$ in the notation. The twoheadarrow \lq $\twoheadrightarrow$\rq will be used to denote the morphisms in $\fib(\XXX)$.\\ 
(2) A functor $f:\XXX \to \YYY$ between categories with fibrations is called {\it{exact}} if it preserves all relevant structures, that is, if $f(\fib(\XXX)) \subset \fib(\YYY)$ and if $f$ preserves pull back along fibration. The last condition means that the canonical map $f(X \times_Z Y) \to f(X) \times_{f(Z)} f(Y)$ is an isomorphism whenever $X \twoheadrightarrow Z$ is in $\fib(\XXX)$. We will denote the category of small category with fibrations and exact functors by $\Fib$.\\
(3) A {\it{category with fibrations and equivalence relations}} $\XXX$ is a category with fibrations whose underlying category $\XXX$ is a category with equivalence relations.\\
(4) A morphism of categories with fibrations and equivalence relations $f:\XXX \to \YYY$ is a functor which is both an exact functor and a morphism of categories with equivalence relations. We will denote the category of small categories with fibrations and equivalence relations by $\Fib_{\er}$.
\end{df}

\begin{ex} \label{cat with fib ex}$ $
\\
(1) An exact category in the sense of Quillen \cite{Qui73} (see also \textbf{Notation \ref{Q-const}} (1)) can be considered as a category with fibrations by choosing a zero object, and defining the fibrations to be the admissible epimorphisms.\\
(2) Let $\XXX$ be a category with fibrations and $f:\YYY \to \XXX$ be a category equivalence. Then $\YYY$ can be considered as a category with fibrations such that $f:\YYY \to \XXX$ is an exact functor by the obvious way.\\
(3) $\AAA$ in {\textbf{Example} \ref{ex of cat with er}} (5) can be considered as a category with fibrations by choosing a zero object, and defining the fibrations to be the surjections.
\begin{proof}
Obviously $\AAA$ is closed under the direct sum operations and taking submodules. Since a fiber product is a submodule of a certain direct sum module, these assertions above implies that $\AAA$ is also closed under taking fiber products. To check other axioms is trivial.
\end{proof} 
\end{ex}

\noindent
We will use the dual $S_{\bullet}$-construction in \cite{Wal85}, denoted by $S^{\op}_{\bullet}$. 

\begin{df}[$S^{\op}_{\bullet}$-construction] \label{S-construction df} $ $
\\
(c.f \cite{Wal85}, \cite{TV04})
Let $\XXX$ be a small category with fibrations.\\
(1) A sequence in $\XXX$
$$X \overset{i}{\to} Y \overset{p}{\twoheadrightarrow} Z$$
is called a {\it{fibration sequence}}, if $p$ is a fibration, and if the following diagram is a cartesian square.
$$\xymatrix{
X \ar[r]^i \ar[d] & Y \ar@{->>}[d]^{p}\\
\ast \ar[r] & Z
}$$
(2) Let $n$ be a non-negative integer. Consider the {\it{arrow category}} $\Ar [n]:=\HOM([1],[n])$ of $[n]$, that is, it is a partially ordered set of pairs
$(i,j)$, $0 \leqq i \leqq j \leqq n$, where $(i,j) \leqq (i',j')$ if and only if $i \leqq i'$ and $j \leqq j'$. We will denote $(i,j)$ as $j/i$.\\
(3) A functor $X:(\Ar[n])^{\op} \to \XXX$ is called the {\it{$n$-th exact functor}} if it satisfies \textbf{Ex 1} and \textbf{Ex 2} below\\
\textbf{Ex 1} $X(i/i)=\ast$ for any $0 \leqq i \leqq n$.\\
\textbf{Ex 2} For any $0 \leqq i \leqq j \leqq k \leqq n$, 
$ X(k/j) \to X(k/i) \to X(j/i)$ is a fibration sequence.\\
(4) A morphism between $n$-th exact functors is a natural transformation. $n$-th exact functors and their morphisms form a category $S_n^{\op}\XXX$. $iS_n^{\op}\XXX$ denotes the subcategory of $S_n^{\op}\XXX$. By definition here a morphism $f:X \to Y$ of $S^{\op}_n \XXX$ is in $iS^{\op}_n \XXX$ if and only if the $f(i/j):X(i/j) \to Y(i/j)$ is an isomorphism for every $0 \leqq j \leqq i \leqq n$. Since $S^{\op}_n\XXX$ is subcategory of $\HOM((\Ar[n])^{\op},\XXX)_{\strict}$, $S_{\bullet}^{\op}\XXX$ and $iS^{\op}_{\bullet}\XXX$ can be considered as sub simplicial categories of $\HOM((\Ar[-])^{\op},\XXX)_{\strict}$. Since this construction is functorial in $\XXX$ by \textbf{Remark \ref{LAX hom is category with equivalence relations}} (4), we have functors
$$ S^{\op}_{\bullet},\ iS^{\op}_{\bullet}:\Fib \to \delcat.$$
(5) We will define a simplicial set $s^{\op}_{\bullet}\XXX$ as $[n] \mapsto \Ob S^{\op}_n\XXX$. Since this construction is functorial in $\XXX$, we have a functor
$$ s^{\op}_{\bullet}:\Fib \to \delset.$$
\end{df}

\begin{rem} \label{functoriality of S-construction} $ $
\\
Let $\XXX$ be a category with fibrations and equivalence relations. Then
since for each non-negative integer $n$, $S^{\op}_n\XXX$ is a subcategory of 
$\HOM((\Ar[n])^{\op},\XXX)_{\strict}$, it can be considered as a category with equivalence relations by \textbf{Example \ref{ex of cat with er}} (2) and \textbf{Remark \ref{LAX hom is category with equivalence relations}} (2). So we get the functors
$$iS^{\op}_{\bullet},\ S^{\op}_{\bullet}:\Fib_{\er} \to \Delta^{\op}\Cat_{\er}.$$ 
\end{rem}

\begin{rem}[Fundamental properties for $s^{\op}_{\bullet}$-construction] \label{fundamental prop of s-const} $ $
\\
Let $f:\XXX \to \YYY$ be an exact functor between small categories with fibrations. Then we have the following properties.\\
(1) The simplicial set $s_{\bullet}\XXX$ is considered as a discrete simplicial category. Then we have the natural inclusion $s_{n}^{\op}\XXX \ni X \mapsto  X \in iS_{n}^{\op}\XXX$ for each $n$. This map is a homotopy equivalence. (See the dual version, \cite{Wal85} p.335 Corollary (2) and p.376, \cite{Sch87} p.272 (1.2) Proposition, see also \cite{HS85} p.416 (3.6).)\\
(2) If $f$ is a category equivalence, that is, it is fully faithful and essentially surjective, then $s^{\op}_{\bullet}f$ is a homotopy equivalence. (See the dual version, \cite{Wal85} p.335 Corollary (1).)\\
\end{rem}

\begin{df}[Cogluing axiom] \label{cogluing axiom} $ $
\\
Let $\XXX$ be a category with fibrations and equivalence relations. We say $\XXX$ satisfies the {\it{cogluing axiom}} if in the following commutative diagrams for $i=1,\ 2$
$$\xymatrix{
X \ar[d]^{a_i} \ar@{->>}[r] & Y \ar[d]^{b_i} & Z \ar[d]^{c_i} \ar[l]\\
X' \ar@{->>}[r] & Y' & Z' \ar[l] & ,
}$$
assertions $a_1 \sim a_2$, $b_1 \sim b_2$ and $c_1 \sim c_2$ imply $a_1 \times_{b_1} c_1 \sim a_2 \times_{b_2} c_2$.\\
\end{df}

\begin{ex} \label{ex of cogluing ax}$ $
\\
(1) Let $\XXX$ be a category with fibrations. We consider it as a category with trivial equivalence relations. Then $\XXX$ satisfies the cogluing axiom.\\
(2) $\AAA$ in {\textbf{Example} \ref{ex of cat with er}} (5) satisfies the cogluing axiom.
\begin{proof}
Consider the diagram in {\textbf{Definition} \ref{cogluing axiom}} for the case of $\AAA$, if $Y$ and $Y'$ are $0$-modules, the result is easy. For the free part of $X \oplus Z$ (resp. $X' \oplus Z'$) is the direct sum of each of the free parts and the torsion part is also. For a general case also follows this special case. For fiber products are submodules of certain direct sum modules and equivalence relations compatible with an inclusion map.\\
\end{proof}
\end{ex}

\begin{prop} \label{Lax group representations} $ $
\\
{\rm (1)} Let $G$ be a group and $\XXX$ be a small category with fibrations and equivalence relations which satisfies the cogluing axiom. Then $\LAX(\ulG,\XXX)_{\strict}$ is a category with fibrations where a morphism $f:X \to Y$ in $\LAX(\ulG,\XXX)_{\strict}$ is a fibration if and only if forgetting the action of $G$, it is in $\fib(\XXX)$.\\ 
{\rm (2)} In the notation above, moreover let $\phi:G \to H$ be a group homomorphism and $\psi:\XXX \to \YYY$ an exact functor between categories with fibrations and equivalence relations satisfying the cogluing axiom. Then the induced functors
$$\LAX(\phi,\XXX)_{\strict}:\LAX(\underline{H},\XXX)_{\strict} \to \LAX(\ulG,\XXX)_{\strict}$$ 
$$\LAX(\ulG,\psi)_{\strict}:\LAX(\ulG,\XXX)_{\strict} \to \LAX(\ulG,\YYY)_{\strict}$$
are exact.\\
\end{prop}

\begin{proof}
In a diagram of lax group representations
$$(X,\rho_X) \twoheadrightarrow (Z,\rho_Z) \leftarrow (Y,\rho_Y),$$
we define a lax action $\rho_{X \times_Z Y}:G \to \Aut (X \times_Z Y)$ by $g \mapsto \rho_X(g) \times_{\rho_Z(g)} \rho_Y(g)$. By the cogluing axiom, we learn that $(X \times_Y Z, \rho_{X \times_Z Y})$ is a lax group representation and a pull back of diagram above. To check other axioms and assertions is trivial.  
\end{proof}

\begin{df}[Lax algebraic $K$-theory] \label{Lax algebraic $K$-theory def} $ $
\\
(1) Let $\XXX$ be a small category with fibrations and equivalence relations. We put a pointed space called as a {\it{lax algebraic $K$-theory space}} associated with $\XXX$ in the following way.
$$\mathbb{K}^{\lax}(\XXX):=\Omega|\Delta iS^{\op}_{\bullet}\XXX|$$
where $\Omega$ is the {\it{loop space functor}}, namely $\Omega=\HHom_{\Top_{\ast}}(S^1,?)$. (See \textbf{Example \ref{cart closed cat ex}} $\text{(2)}_{\ast}$).\\
(2) In the notation above, assume $\XXX$ satisfies the cogluing axiom. For any group $G$, we put a pointed space called as a {\it{lax representations space}} associated with $G$ and $\XXX$ in the following way.
$$\mathbb{R}^{\lax}(G,\XXX):=\Omega \B iS^{\op}_{\bullet}\LAX(\ulG,\XXX)_{\strict}$$
where $B$ is the {\it{classifying space functor}}, that is, composition of 
$\Cat \overset{\N}{\to} \delset \overset{|?|}{\to} \Top$. We put its $n$-th homotopy group for any non-negative integer $n$ as
$$R^{\lax}_n(G,\XXX):=\pi_n(\mathbb{R}^{\lax}(G,\XXX))$$
and call it as the {\it{$n$-th lax representations group}} associated with $G$ and $\XXX$.\\
(3) In the notation above, by \textbf{Proposition} \ref{Lax group representations} (2), for any non-negative integer $n$, the canonical homomorphism $G \to \pt$ induces the group homomorphism 
$$R_n(G,\XXX) \to R_n(\pt,\XXX)$$
where $\pt$ means the unit group. We put
$$\tilde{R}_n^{\lax}(G,\XXX):=\Ker(R_n(G,\XXX) \to R_n(\pt,\XXX))$$
and call it as the {\it{$n$-th reduced lax group representation group}} associated with $G$ and $\XXX$.
\end{df}

\begin{exc} $ $
\\
We can define a concept of Waldhausen categories with equivalence relations and associate them with lax algebraic $K$-theory in the obvious way. But we will not need them to prove the main theorem and the author abandon developing these theories and full functoriality of lax algebraic $K$-theory.
\end{exc}

\begin{ex} $ $
\\
In the notation above, we assume $\XXX$ is a category with trivial equivalence relations. Then\\
(1) $\mathbb{K}^{\lax}(\XXX)$ is just a usual algebraic $K$-space associated with $\XXX$. We will denote it by $\mathbb{K}(\XXX)$.\\
(2) For any group $G$ (and non-negative integer $n$), $\mathbb{R}^{lax}(G,\XXX)$ (resp. $R^{\lax}_n(G,\XXX)$, $\tilde{R}^{\lax}_n(G,\XXX)$) is just a usual representation space (resp, $n$-th representation group, $n$-th reduced representation group) associated with $G$ and $\XXX$. We will denote it by $\mathbb{R}(G,\XXX)$ (resp. $R_n(G,\XXX)$, $\tilde{R}_n(G,\XXX)$).
\end{ex}

\begin{rem}[Edgewise subdivison] \label{edgewise subdivision} $ $
\\
We review the {\it{edgewise subdivison}} functor. (c.f. \cite{Seg73} p.219, \cite{Wal85} p.375.)\\
(1) First we define the {\it{doubling functor}} $d:\Delta \to \Delta$ as follows. For each $[n]$, we put $d([n]):=[2n+1]$ and for each ordered map $a:[n] \to [m]$, we put
$$d(a)(k):=
\begin{cases}
m-a(n-k) & \text{if $0 \leqq k \leqq n$}\\
a(k-n+1)+m+1 & \text{if $n+1 \leqq k \leqq 2n+1$}.
\end{cases}$$
(2) For any simplicial object $X:\Delta^{\op} \to \XXX$ in a category $\XXX$, associate another $X^{\edge}:\Delta^{\op} \to \XXX$, namely the composite $X^{\edge}:=Xd^{\op}$. This assignment define the {\it{edgewise subdivision}} functor
$$?^{\edge}:\Delta^{\op}\XXX \to \Delta^{\op}\XXX.$$
(3) If $X$ is a simplicial space then the geometric realizations $|X|$ and $|X^{\edge}|$ are naturally homeomorphic. (c.f. \cite{Seg73} Ibid Proposition (A.1).)\\
(4) For any small category with fibrations and equivalence relations $\XXX$, applying the fact in (3) to the simplicial space $[n] \mapsto |\Delta iS_n^{\op}\XXX|$ we obtain that $\mathbb{K}^{\lax}(\XXX)$ and $\mathbb{K}^{\lax,\edge}(\XXX):=\Omega|\Delta iS_{\bullet}^{\op,\edge}\XXX|$ are homotopy equivalent. (See \cite{Wal85} Ibid.) If $\XXX$ satisfies the cogluing axiom and let $G$ be a group, then we will also use $\mathbb{R}^{\lax,\edge}(G,\XXX):=\Omega BiS^{\op,\edge}_{\bullet}\LAX(\ulG,\XXX)_{\strict}$ which is homotopy equivalent to $\mathbb{R}^{\lax}(G,\XXX)$.
\end{rem}

\begin{thm} \label{rectification variant} $ $
\\
Let $\XXX$ be a category with fibrations and equivalence relations satisfying the cogluing axiom. Then the canonical map induced from the identity functor
$$\mathbb{K}(\XXX) \to \mathbb{K}^{\lax}(\XXX)$$
has a retraction up to homotopy.
\end{thm}

\begin{proof}
First notice the following commutative diagram of trisimplicial sets
$$\xymatrix{
s^{\op}_{\bullet}\XXX \ar[r] & \N iS^{\op}_{\bullet}\XXX \ar[r] \ar[d] & \N\N^{\er}iS^{\op}_{\bullet}\XXX \ar[d]\\
& \Delta iS^{\op}_{\bullet}\XXX \ar[r] & \N\underline{\Delta} iS^{\op}_{\bullet}\XXX 
}$$
where $s^{\op}_{\bullet}\XXX$ is considered as a constant bisimplicial simplicial set (namely trisimplicial set) and so on. By \textbf{Remark \ref{FPVN}} (3), \textbf{Remark \ref{fundamental prop of s-const}} (1), we shall only prove that the inclusion map
$s^{\op}_{\bullet}\XXX \hookrightarrow \N\N^{\er}iS^{\op}_{\bullet}\XXX$ has a retraction. Now we consider a simplicial category $M_{\bullet}\XXX$ which is a sub simplicial category of
$[n] \mapsto \LAX([n] \times [n],\XXX)$ as follows.\\
Object: A lax functor $X:[n] \times [n] \to \XXX$ satisfying the following conditions.\\
(i) For any $0 \leqq i \leqq n$ and $0 \leqq s \leqq t \leqq n$, $X(i,s)=X(i,t)$ and $X(i,s \leqq t)=\id:X(i,s) \to X(i,t)$ and,\\
(ii) for any $0 \leqq j \leqq n$ and $0 \leqq s \leqq n-1$, $X_{s,j}:=X(s \leqq s+1,j):X(s,j) \to X(s+1,j)$ is an isomorphism.\\
Morphism: A deformation $f:X \to Y$ satisfying the following conditions.\\
(iii) For any $0 \leqq i \leqq n$ and $0 \leqq s \leqq t \leqq n$, $f(i,s)=f(i,t)$ and,\\
(iv) for any $0 \leqq j \leqq n$ and $0 \leqq s \leqq n-1$, the following diagram is commutative:
$$\xymatrix{
X(s,j) \ar[d]_{f(s,j)} \ar[r]^{X_{s,j}}_{\!\!\!\!\!\sim} & X(s+1,j) \ar[d]^{f(s+1,j)}\\
Y(s,j) \ar[r]_{Y_{s,j}}^{\!\!\!\!\!\sim} & Y(s+1,j). 
}$$
\textbf{Claim}\\
(1) For each $n$, $M_n\XXX$ can be considered as a category with fibrations where a morphism is a fibration if and only if it is a term-wised fibrations in $\XXX$. \\
(2) For each ordered map $\phi:[n] \to [m]$, the induced map $M_m\XXX \to M_n\XXX$ is exact.\\
(3) $\diag\N\N^{\er}iS^{\op}_{\bullet}\XXX=\diag s^{\op}_{\bullet}M_{\bullet}\XXX$. Here $\diag$ means the {\it{diagonalization functor}}.\\
\begin{proof}[Proof of \textbf{Claim}]
To prove \textbf{Fib 3}, we need the assumption that $\XXX$ satisfies the cogluing axiom. Fiber product of a pull back diagram in $M_n\XXX$ is a term-wised pull back. To prove other axioms and assertions is trivial.
\end{proof} 
\noindent
Let $M'_{\bullet}\XXX$ be a full sub simplicial category of $M_{\bullet}\XXX$ whose $n$-simplexes objects are lax functors $X:[n]\times[n] \to \XXX$ such that $X(i,j)=X(s,t)$ for any $0 \leqq i,j,s,t, \leqq n$.\\
\\
\noindent
\textbf{Claim 2}\\
For each non-negative integer $n$, the inclusion functor $M'_n\XXX \hookrightarrow M_n\XXX$ is essentially surjective.
Hence by \textbf{Example \ref{cat with fib ex}} (2), $M'_{\bullet}\XXX$ is also a simplicial category with fibrations.\\
\begin{proof}[Proof of \textbf{Claim 2}]
For any object $X$ in $M_n\XXX$, we define an object $Z$ in $M'_n\XXX$ and an isomorphism $f:Z \isoto X$ as follows.\\
First we put $Z(i,j):=X(0,0)$ and
$Z_{p,q}:=X_{0,0}^{-1}X_{1,0}^{-1}\cdots X_{p,0}^{-1}X_{p,q}X_{p-1,0}X_{p-2,0}\cdots X_{0,0}$ for any $0 \leqq i,j \leqq n$ and $0 \leqq p,q \leqq n-1$ where $X_{k,0}$ for $k<0$ means $\id$. We can easily check $Z$ is an object in $M'_n\XXX$. Next we put $f(i,j):=X_{i-1,0}X_{i-2,0}\cdots X_{0,0}$ for any $0 \leqq i,j \leqq n$. We can easily check that this define an isomorphism $f:Z \isoto X$.
\end{proof} 
So we have a homotopy equivalence $s^{\op}_{\bullet}M'_{\bullet}\XXX \to s^{\op}_{\bullet}M_{\bullet}\XXX$ by \textbf{Remark \ref{fundamental prop of s-const}} (2) and the \textbf{Realization lemma \ref{real lem}} below. Since the canonical map factor through $s^{\op}_{\bullet}\XXX \to s^{\op}_{\bullet}M'_{\bullet}\XXX \to s^{\op}M_{\bullet}\XXX$, we shall only prove that $s^{\op}_{\bullet}\XXX \to s^{\op}_{\bullet}M'_{\bullet}\XXX$ has a retraction.\\
\\
\noindent
\textbf{Claim 3}\\
$s^{\op}_{\bullet}\XXX \to s^{\op}_{\bullet}M'_{\bullet}\XXX$ has a retraction.\\
\begin{proof}[Proof of \textbf{Claim 3}]
For any non-negative integer $n$, an element in $(\diag \N M'_{\bullet}\XXX)_n$ is consist of a family $\{X,a_{i,j}\}_{0 \leqq i,j \leqq n}$ where $X$ is the $n$-th exact functor $(\Ar[n])^{\op} \to \XXX$ and $a_{i,j}$ are automorphisms of $X$ such that $a_{i,p} \sim a_{i,q}$ for any $0 \leqq i, p, q \leqq n$. So we have an assignment 
$$(\diag \N M'_{\bullet}\XXX)_n \ni \{X,a_{i,j}\}_{0 \leqq i,j \leqq n} \mapsto X \in s^{\op}_n\XXX$$  
which defines a desired retraction.
\end{proof}
Hence the canonical map $\mathbf{K}(\XXX) \to \mathbf{K}^{\lax}(\XXX)$ has a retraction up to homotopy.
\end{proof}

\begin{lem} [Realization lemma] \label{real lem} $ $
\\
{\rm (c.f. \cite{Seg74} Appendix A, \cite{Wal78} p. 164 Lemma 5.1).}\\
Let $X_{\bullet \bullet} \to Y_{\bullet \bullet}$ be a map of bisimplicial sets. Suppose that for every $n$, the map $X_{\bullet n} \to Y_{\bullet n}$ is a homotopy equivalence. Then $X_{\bullet \bullet} \to Y_{\bullet \bullet}$ is a homotopy equivalence. 
\end{lem}

\section{The lax Sherman map}

In this subsection, we will construct the {\it{lax Sherman map}} and by using them, state the universal property of algebraic $K$-theory associated with semi-simple exact categories.

\subsection{The exponential law}

As explained in \S 1, in this section, we will study behaviour of internal hom objects through the equivalences in (\ref{dogma}). To do so, we will start to define {\it{cartesian closed categories}}.  

\begin{df}[Cartesian closed categories] $ $
\\
Let $\VVV$ be a category with a terminal object $\ast$, and finite products. We pick a functor 
$$\times:\VVV \times \VVV \to \VVV$$ 
sending $V_1,V_2$ to $V_1 \times V_2$, the product of $V_1$ and $V_2$. Suppose that $\VVV$ provided the internal hom functor, i.e., a functor 
$\HHom_{\VVV}:\VVV^{op} \times \VVV \to \VVV$ and we have an isomorphism, natural in $U,V$ and $W \in \VVV$ of morphism sets:
$$\exp:\Hom(U \times V,W) \overset{\sim}{\to} \Hom(U,\HHom_{\VVV}(V,W)).$$     
We call such a $\VVV$ together with the structure $\times$, $\HHom_{\VVV}(\ ,\ )$, $\exp$, a {\it{cartesian closed category}}.   
\end{df}

\begin{ex}  \label{cart closed cat ex} $ $
\\
(1) $\VVV=\Set$, the category of sets, where $\HHom_{\VVV}(U,V)$ is the set of functions $U \to V$.\\
(2) $\VVV=\Top$, the category of compactly generated spaces, where $\HHom_{\VVV}(U,V)$ is the space of continuous functions $U \to V$, with the compact-open topology strengthened to a compactly generated topology.\\
$\text{(2)}_{\ast}$ $\VVV=\Top_{\ast}$, the category of pointed compactly generated spaces with non-degenerated base points, where $\HHom_{\VVV}(U,V)$ is the subspace of $\HHom_{\Top}(U,V)$ whose maps are base point preserving.\\
(3) $\VVV=\Delta^{\op} \Set$, the category of simplicial sets, where $\HHom_{\VVV}(U,V)$ is the usual simplicial function space
$[n] \mapsto \Hom_{\delset}(N[n] \times U,V)$.\\
(4) $\VVV=\Cat$, the category of small categories, where $\HHom_{\VVV}(U,V)$ is the category with objects functors $U \to V$, and morphisms the natural transforms, that is, $\HOM(U,V)_{\strict}$\\
\end{ex}

\begin{lem} \label{preserveness of int hom} $ $
\\
Let $f:\VVV \to \WWW$ be a functor between cartesian closed categories. Assume $f$ preserves products. Then we have a canonical natural transformation
$$f(\HHom_{\VVV}(?,?)) \Rightarrow \HHom_{\WWW}(f(?),f(?))$$ 
between the functors $\VVV^{\op} \times \VVV \to \WWW$.
\end{lem}

\begin{proof}
For each objects $V_1$, $V_2$ in $\VVV$, we have the {\it{evaluation morphism}}
$$\ev_{V_1,V_2}:\HHom_{\VVV}(V_1,V_2) \times V_1 \to V_1$$
which is functorial in $V_1$ and $V_2$ and corresponded to $\id_{\HHom(V_1,V_2)}$ through the following isomorphism
$$\exp:\Hom_{\VVV}(\HHom_{\VVV}(V_1,V_2) \times V_1,V_2) \isoto \Hom_{\VVV}(\HHom_{\VVV}(V_1,V_2),\HHom_{\VVV}(V_1,V_2)).$$
Since $f$ preserves products, we have the following isomorphisms
\begin{multline*}
\Hom_{\WWW}(f(\HHom_{\VVV}(V_1,V_2) \times V_1),f(V_2)) \isoto
\Hom_{\WWW}(f(\HHom_{\VVV}(V_1,V_2)) \times f(V_1),f(V_2))\\
 \overset{\exp}{\isoto}\Hom_{\WWW}(f(\HHom_{\VVV}(V_1,V_2)),\HHom_{\WWW}(f(V_1),f(V_2)).
\end{multline*}
Hence we get the morphism 
$$f(\HHom_{\VVV}(V_1,V_2)) \to \HHom_{\WWW}(f(V_1),f(V_2))$$
which is corresponded to $f(\ev_{V_1,V_2})$ through the isomorphisms above, functorial in $V_1$ and $V_2$ and the desired natural transformation. 
\end{proof}

\begin{ex} \label{Sherman argument} $ $
\\
(1) It is well-known that the nerve functor $\N:\Cat \to \delset$ preserves products. So by \textbf{Lemma \ref{preserveness of int hom}}, we have the canonical natural transformation
$$\N\HOM(?,?)_{\strict} \Rightarrow \HHom_{\delset}(\N(?),\N(?)).$$ 
(Compare the argument in \cite{She82} p.212.) We will construct this \lq lax version\rq \ in \textbf{Example \ref{lax exp law cor 3}} (1).\\
(2) It is well-know that the geometric realization functor $|?|:\delset \to \Top$ preserves products. (See \cite{Mil57}.) So by \textbf{Lemma \ref{preserveness of int hom}}, we have the canonical natural transformation
$$|\HHom_{\delset}(?,?)| \Rightarrow \HHom_{\Top}(|?|,|?|).$$
(Compare the argument in Ibid.) We will construct this \lq bisimplicial version\rq \ in \textbf{Example \ref{lem example}}.
\end{ex}

\noindent
The following \textbf{Lemma \ref{internal hom lem}} can be easily verified.

\begin{lem} \label{internal hom lem}$ $
\\
Let $I$ be a category, $\VVV$ a cartesian closed category and $f:I \to \VVV$ a functor. For each object $V$ in $I$, we consider it as a constant functor
$V:I \to \VVV$. Then the following functor
$$\HHom(V,f):I \ni i \mapsto \HHom_{\VVV}(V,f(i)) \in \VVV$$
is an internal hom object in $\HOM(I,\VVV)$.
\end{lem}

\begin{ex} \label{lem example} $ $
\\
Let $X$ be a simplicial set and $Y$ a bisimplicial set. Then we have a simplicial simplicial set $[n] \mapsto \HHom_{\delset}(X,Y_{n \bullet})$. We define a continuous map
$$|\diag \HHom_{\delset}(X,Y)| \to \HHom_{\Top}(|X|,|Y|)$$
natural in each variable, as follows.\\
First notice that we have two simplicial topological spaces $[n] \mapsto |\HHom_{\delset}(X,Y_{n \bullet})|$ and $[n] \mapsto \HHom_{\Top}(|X|,|Y_{n \bullet}|)$. By \textbf{Example \ref{Sherman argument}} (2), we have the natural transformation between them. Therefore we also have that between their geometric realizations
$$|[n] \mapsto |\HHom_{\delset}(X,Y_{n \bullet})|| \to |[n] \mapsto \HHom_{\Top}(|X|,|Y_{n \bullet}|)|.$$
Notice that by \textbf{Lemma \ref{internal hom lem}}, $[n] \mapsto \HHom_{\Top}(|X|,|Y_{n \bullet}|)$ is an internal hom object in $\Delta^{\op}\Top$ and it is well-known that $|?|:\Delta^{\op}\Top \to \Top$ also preserves products. Hence we have the following map 
$$|[n] \mapsto \HHom_{\Top}(|X|,|Y_{n \bullet}|)| \to \HHom_{\Top}(|X|,|[n] \mapsto |Y_{n \bullet}|)$$
by \textbf{Lemma \ref{preserveness of int hom}}. Composing the two maps above, we get the desired map. Notice that if $Y$ is a constant simplicial simplicial set, then the construction above is coincided with that in \textbf{Example \ref{Sherman argument}} (2).    
\end{ex}

\begin{thm}[Lax exponential law] \label{Lax exponential law} $ $
\\
Let $I$, $J$ be small categories and $\XXX$ a category with equivalence relations. Then we have an \lq almost\rq \ isomorphism as categories with equivalence relations
$$\LAX(I,\LAX(J,\XXX)) \isoto \LAX(I \times J,\XXX)$$
which is functorial in each variable.
\end{thm}

\begin{proof}
The phrase \lq almost\rq \ in the statement means the following sense.\\
There are two morphisms of categories of equivalence relations
$$A:\LAX(I,\LAX(J,\XXX)) \to \LAX(I \times J,\XXX),$$
$$B:\LAX(I \times J,\XXX) \to \LAX(I,\LAX(J,\XXX))$$
such that $BA=\id$ and there is a deformation $\Theta:AB \Rightarrow \id$ relative to $\Ob \LAX(I \times J,\XXX)$.\\
The definition of $A$:\\
For any lax functors $f,\ g:I \to \LAX(J,\XXX)$, any deformation $x: f \Rightarrow g$ and any morphism $(a,b):(i_1,j_1) \to (i_2,j_2)$ in $I \times J$, we put 
$$\begin{cases}
A(f)(i_1,j_1):=f(i_1)(j_1)\\
A(f)(a,b):=f(a)(j_2)f(i_1)(b)\\
A(x)(i_1,j_1):=x(i_1)(j_1). 
\end{cases}$$
(In the definition above, we may replace $A(f)(a,b)=f(a)(j_2)f(i_1)(b)$ with
$A(f)(a,b):=f(i_2)(b)f(a)(j_1)$.)\\
The definition of $B$:\\
$f,\ g:I \times J \to \XXX$, any deformation $x: f \Rightarrow g$ and any morphism $a:i_1 \to i_2$ in $I$ and $j_1 \to j_2$ in $J$, we put 
$$\begin{cases}
B(f)(i_1)(j_1):=f(i_1,j_1)\\
B(f)(a)(j_1):=f(a,j_1)\\
B(f)(i_1)(b):=f(\id_i,b)\\
B(x)(i_1,j_1):=x(i_1,j_1). 
\end{cases}$$
The definition $\Theta$ is obvious. Then we can easily verify the assertions above.
\end{proof}

\begin{cor} \label{Lax exponential law cor}$ $
\\
As in the notation \textbf{Theorem \ref{Lax exponential law}}, \lq almost\rq \ isomorphisms
$$\LAX(I,\LAX(J,\XXX)) \isoto \LAX(I \times J,\XXX)
\isoto \LAX(J,\LAX(I,\XXX))$$
induce an actual isomorphism
$$\HOM(I,\LAX(J,\XXX)_{\strict})_{\strict} \isoto  \LAX(J,\HOM(I,\XXX)_{\strict})_{\strict}$$
natural in each variable.
\end{cor}

\begin{proof}
The proof is straightforward. Notice that in the proof of \textbf{Theorem \ref{Lax exponential law}}, we have two choice of defining desired isomorphisms. But these isomorphisms induce same isomorphism on $\HOM(I,\LAX(J,\XXX)_{\strict})_{\strict}$.
\end{proof}

\begin{cor} \label{lax exp law cor 2} $ $
\\
Let $G$ be a group and $\XXX$ a category with fibrations and equivalence relations satisfying the cogluing axiom. Then by \textbf{Proposition} \ref{Lax group representations}, $\LAX(\ulG,\XXX)$ is a category with fibrations and equivalence relations. We have a canonical isomorphism of simplicial categories
$$iS^{\op}_{\bullet}\LAX(\ulG,\XXX)_{\strict} \isoto \LAX(\ulG,iS^{\op}_{\bullet}\XXX)_{\strict}$$
natural in each variable.
\end{cor}

\begin{proof}
For each non-negative integer $n$, the isomorphism
$$\HOM((\Ar[n])^{\op},\LAX(\ulG,\XXX)_{\strict})_{\strict} \isoto  \LAX(\ulG,\HOM((\Ar[n])^{\op},\XXX)_{\strict})_{\strict}$$
induces an isomorphism
$$iS^{\op}_{n}\LAX(\ulG,\XXX)_{\strict} \isoto \LAX(\ulG,iS^{\op}_{n}\XXX)_{\strict}$$
which is functorial in $n$ and each variable and this is the desired isomorphism.
\end{proof}

\begin{ex} \label{lax exp law cor 3}  $ $
\\
(1) Let $I$ be a small category and $\XXX$ a category with equivalence relations. We have a simplicial map
$$\N\LAX(I,\XXX)_{\strict} \to \HHom_{\delset}(\N I,\Delta \XXX)$$ 
natural in each variable, as follows.\\
For each non-negative integer $n$, by \textbf{Theorem \ref{Lax exponential law}}, we have the canonical morphism
\begin{multline*}
\N_n\LAX(I,\XXX)_{\strict}=\Ob \HOM([n],\LAX(I,\XXX)_{\strict})=\Hom_{\Cat}([n],\LAX(I,\XXX)_{\strict})\\
\to \Ob\LAX([n] \times I,\XXX)=\Lax([n] \times I,\XXX)
\end{multline*}
and by \textbf{Remark \ref{FPVN}} (2), we have the canonical morphism
$$\Lax([n] \times I,\XXX) \to \Hom_{\delset}(\Delta [n] \times \N I, \Delta \XXX).$$
Composing these morphisms, we get the desired morphism. Notice that if $\XXX$ is a trivial category with equivalence relations, then the construction above is coincided with that in \textbf{Example \ref{Sherman argument}} (1).\\
(2) Let $G$ be a group and $\XXX$ a small category with fibrations and equivalence relations satisfying the cogluing axiom. Then composing the maps in (1) and \textbf{Corollary \ref{lax exp law cor 2}}, we get the following morphism
$$\N iS^{\op}_{\bullet}\LAX(\underline{G},\XXX)_{\strict} \isoto \N\LAX(\ulG,iS^{\op}_{\bullet}\XXX)_{\strict} \to
\HHom_{\delset}(\N\ulG,\Delta iS^{\op}_{\bullet}\XXX)$$
natural in each variable.
\end{ex}

\begin{conc} \label{pre lax sherman} $ $
\\
Let $G$ be a group and $\XXX$ be a small category with fibrations and equivalence relations satisfying the cogluing axiom. Then composing morphisms 
\begin{multline*}
\mathbb{R}^{\lax,\edge}(G,\XXX)=\Omega BiS^{\op,\edge}_{\bullet}\LAX(\ulG,\XXX)_{\strict} \to \Omega |\HHom_{\delset}(\N\ulG,\Delta iS^{\op,\edge}_{\bullet}\XXX)|\\
\to \Omega \HHom_{\Top}(BG,|\Delta iS^{\op,\edge}_{\bullet}\XXX|)
\end{multline*}
in \textbf{Example \ref{lem example}} and \textbf{Example \ref{lax exp law cor 3}}, we get the morphism
$$\mathbb{R}^{\lax,\edge}(G,\XXX) \to \Omega \HHom_{\Top}(BG,|\Delta iS^{\op,\edge}_{\bullet}\XXX|)$$
natural in each variable. 
\end{conc}

\noindent
Before constructing the {\it{lax Sherman map}}, we shall review about topological spaces. So we will do in the next subsection. 

\subsection{The construction}

In this subsection, we will define the {\it{lax Sherman map}}. To do so, we will start to review fundamental facts about topologies. 

\begin{nt} \label{Top not} $ $
\\
(1) The set of free homotopy classes of free maps from $X$ to $Y$ will be denoted $[X,Y](=\pi_0(\HHom_{\Top}(X,Y)))$, while $[X,Y]_{\ast}(=\pi_0(\HHom_{\Top_{\ast}}(X,Y)))$ will denote the set of pointed homotopy classes of pointed maps.\\
(2) For a pointed space $X$, $X_0$ will denote the path-component containing the basepoint.\\
\end{nt}

\begin{rev}[Well-known result for function spaces] $ $
\\
(1) If $Y$ is a connected simple space, then the canonical map $[X,Y]_{\ast} \to [X,Y]$ is a bijection. In particular, in this case, if two pointed maps $f_1,f_2:X \to Y$ are freely homotopic, then they are homotopic by a pointed homotopy.\\
(2) If $Y$ is a connected pointed space, then we have a canonical bijection
$\HHom_{\Top_{\ast}}(Y,X_0) \to \HHom_{\Top_{\ast}}(Y,X)$.\\
\end{rev}

\begin{df}[$H$-spaces] \label{H-sp def} $ $
\\
(1) An {\it{$H$-spaces}} will be a pointed space $X$ equipped with a pointed addition map $\mu:X \times X \to X$ which, up to homotopy, is associative and has the base point as unit element.\\
(2) An {\it{$H$-homomorphism}} of $H$-spaces will be pointed map preserving the addition up to (pointed) homotopy.\\
(3) An {\it{$H$-isomorphism}} will be an $H$-homomorphism which is a homotopy equivalence.\\
(4) A {\it{homotopy inverse}} for an $H$-space $X$ will be a pointed map $X \to X$ which is an inverse for the addition, up to pointed homotopy.
\end{df}

\begin{rev}[Well-known results for $H$-spaces]  \label{H-space fact}$ $
\\
Let $X$ be an $H$-space.\\
(1) If $X$ has a homotopy inverse, then $\pi_0(X)$ is a group. Conversely, if $X$ is a CW-complex and $\pi_0(X)$ is a group, then $X$ has a homotopy inverse.\\ 
(2) The addition restricts to an addition on $X_0$. If $X$ has a homotopy inverse, then it restricts to a homotopy inverse for $X_0$.\\
Moreover let $Y$ be an arbitrary pointed space.\\
(3) Pointwise addition of functions defines $H$-space structures on both $\HHom_{\Top}(Y,X)$ and $\HHom_{\Top_{\ast}}(Y,X)$. Moreover, if $X$ has a homotopy inverse, then so does each of these function spaces. In this case, $[Y,X]$ and $[Y,X]_{\ast}$ are groups.
\end{rev}

\begin{rev}[$2$-coskeleton map]  \label{2-cosk def} $ $
\\
Let $X$ be a connected pointed CW-complex and $Y$ be a $H$-spaces.\\
(1) There is a canonical pointed homotopy class of maps $X \to B\pi_1(X)$ (the {\it{$2$-coskeleton}}); more precisely, this is the class corresponding to the identity map of $\pi_1(X)$ under the isomorphism 
$$[X,B\pi_1(X)]_{\ast} \isoto \Hom(\pi_1(X),\pi_1(X)).$$
(2) Composition with this define a homomorphism
$$[B\pi_1(X),Y] \to [X,Y].$$
\end{rev}

\begin{df}[Lax Sherman map] \label{lax Sherman df}$ $
\\
(1) In the notation \textbf{Conclusion \ref{pre lax sherman}}, taking homotopy group we get a homomorphism, for any non negative integer $n$:
\begin{eqnarray*}
R_n^{\lax}(G,\XXX) &=& \pi_{n}(\mathbb{R}^{\lax,\edge}(G,\XXX))\\ 
& \to & \pi_{n}(\Omega \HHom_{\Top}(BG,|\Delta iS_{\bullet}^{\op,\edge}\XXX|))\\
& \isoto & \pi_0(\Omega^{n+1}\HHom_{\Top}(BG,|\Delta iS_{\bullet}^{\op,\edge}\XXX|))\\
& \isoto & \pi_0(\HHom_{\Top}(BG,\Omega^{n}\mathbb{K}^{\lax,\edge}(\XXX)))\\
& = & [BG,\Omega^{n}\mathbb{K}^{\lax,\edge}(\XXX)].
\end{eqnarray*}
(Here we are using the fact that the canonical homomorphism
$$\Omega^{n+1}\HHom_{\Top}(BG,|\Delta iS_{\bullet}^{\op,\edge}\XXX|)\isoto\HHom_{\Top}(BG,\Omega^{n+1}|\Delta iS_{\bullet}^{\op,\edge}\XXX|)$$ 
is an $H$-homomorphism, hence induces an isomorphism on $\pi_0$.)\\
This map is clearly natural in each variable.\\
(3) For any pointed connected CW-complex $X$, putting this together with the $2$-coskeleton map above, we get a homomorphism
$$\Sh^{\lax}(X,\XXX)_n:R_n^{\lax}(\pi_1(X),\XXX) \to [X,\Omega^{n}\mathbb{K}^{\lax,\edge}(\XXX)],$$
representing a natural transformation of bifunctors 
$$\HHH(\mathbf{CW}_{\ast}) \times \Fib_{\er} \to \mathbf{Grp}$$
where $\HHH(\text{\textbf{CW}}_{\ast})$ is the pointed homotopy category of pointed CW-complexes and $\mathbf{Grp}$ is the category of groups.
\end{df}

\begin{var}[Reduced cases] \label{variant}$ $
\\
Let $\pt$ denote the trivial space.\\
(1) Since $B\pi_1(X)$ is connected, the kernel of the homomorphism
$$[B\pi_1(X),\Omega^{n}\mathbb{K}^{\lax,\edge}\XXX] \to [\pt,\Omega^{n}\mathbb{K}^{\lax,\edge}\XXX] \isoto \pi_0(\Omega^{n}\mathbb{K}^{\lax,\edge}\XXX)$$
is clearly $[B\pi_1(X),{(\Omega^{n}\mathbb{K}^{\lax,\edge}\XXX)}_0]$.\\
(3) Furthermore, since ${(\Omega^{n}\mathbb{K}^{\lax,\edge}\XXX)}_0$ is an $H$-space (hence simple) and connected, the canonical map 
$$[B\pi_1(X),{(\Omega^{n}\mathbb{K}^{\lax,\edge}\XXX)}_0]_{\ast} \to [B\pi_1(X),{(\Omega^{n}\mathbb{K}^{\lax,\edge}\XXX)}_0]$$ 
is an isomorphism.\\
(4) It follows that there is an induced map
$$\tilde{R}_n^{\lax}(\pi_1(X),\XXX) \to [B\pi_1(X),{(\Omega^{n}\mathbb{K}^{\lax,\edge}\XXX)}_0]_{\ast},$$
and thus also a map
$$\tilde{\Sh}^{\lax}(X,\XXX)_n:\tilde{R}_n^{\lax}(\pi_1(X),\XXX) \to [X,{(\Omega^{n}\mathbb{K}^{\lax,\edge}\XXX)}_0]_{\ast}.$$
\end{var}

\begin{rem} $ $
\\
In the notation above, if equivalence relations in $\XXX$ are trivial, we get maps
$$\Sh(X,\XXX)_n:R_n(G,\XXX) \to [BG,\Omega^n\mathbb{K}^{\edge}(\XXX)],$$
$$\tilde{\Sh}(X,\XXX)_n:\tilde{R}_n(G,\XXX) \to [BG,{(\Omega^n\mathbb{K}^{\edge}(\XXX))}_0]_{\ast}.$$
\end{rem}

\subsection{Compatibility with the original one and the universal property}

In this subsection, we will prove that the lax Sherman map is compatible with classical one defined in \cite{She82} in the suitable case. First we will review the dual notion of $Q$-construction in \cite{Qui73}.

\begin{nt}[$Q^{\op}$-construction] \label{Q-const} $ $
\\
(1) First we shall quickly review the definition of {\it{exact categories}}. An exact category is a triple of $(\EEE,\fib(\EEE),\cof(\EEE))$ consisting of an additive category $\EEE$ with  subcategories $\fib(\EEE),\ \cof(\EEE) \hookrightarrow \EEE$ whose morphisms will be called {\it{admissible epimorphisms}} and {\it{admissible monomorphisms}} respectively satisfying the following axioms:\\
(a) Pairs $(\EEE,\fib(\EEE))$ and $(\EEE^{\op},(\cof(\EEE))^{\op})$ are categories with fibrations.\\
(b) A sequence in $\EEE$
$$X\overset{i}{\to}Y\overset{p}{\to}Z$$   
is a fibration sequence in $(\EEE,\fib(\EEE))$ if and only if it is in $(\EEE^{\op},(\cof(\EEE))^{\op})$.\\
(Compare \cite{Kel90} p.405 A.1, this definition is equivalent to Quillen's one \cite{Qui73}. See also \cite{TT90} p.253 1.2.2). We will always omit $\fib(\EEE)$ and $\cof(\EEE)$ in the notation.\\
(2) For each small exact category $\EEE$, we form a new category $Q^{\op}\EEE$, the dual notion in \cite{Qui73}. That is, $Q^{\op}\EEE$ has the same objects as $\EEE$ and morphisms are defined in the following way. A morphism from $X$ to $Y$ in the category $Q^{\op}\EEE$ is by definition an isomorphism class of the following diagrams
$$ X \overset{i}{\rightarrowtail} W \overset{p}{\twoheadleftarrow} Y$$
where $i$ is an admissible monomorphism and $p$ is an admissible epimorphism. Here we consider isomorphisms of these diagrams which induces the identity on $X$ and $Y$.\\
Given a morphism from $Y$ to $Z$ represented by the diagram
$$Y \overset{j}{\rightarrowtail} V \overset{q}{\twoheadleftarrow} Z$$
the composition of this morphism with the morphism from $X$ to $Y$ above is the morphism represented by the pair $\iota_1 i$, $\iota_2 q $ in the diagram:
$$\xymatrix{
W \coprod_{Y} V & V \ar[l]_{\ \ \ \iota_2} & Z \ar[l]_{q}\\
W \ar[u]^{\iota_1} & Y \ar[l]_{\ \ \ p} \ar[u]_{j}\\
X \ar[u]^{i} & &.
}$$
It is clear that composition is well-defined and associative.\\
(3) In the notation above, we also form a new simplicial category $iQ^{\op}_{\bullet} \EEE$, the dual notion in \cite{Wal85} p.375. $iQ^{\op}_{\bullet} \EEE$ is a sub simplicial category of $[n] \mapsto \HOM([n],Q^{\op} \EEE)_{\strict}$ whose morphisms are all isomorphisms. 
\end{nt}

\noindent
We will review the original definion of the Sherman map for exact categories.

\begin{ex}[The original Sherman map] \label{origin} $ $
\\
We will review the original construction of the Sherman maps. Let fix an small exact category $\EEE$.\\
(1) For any group $G$, as a variant of \textbf{Corollary \ref{lax exp law cor 2}}, we have the following canonical identities
$$iQ^{\op}_{\bullet}\HOM(\ulG,\EEE)_{\strict} \isoto \HOM(\ulG,iQ^{\op}_{\bullet}\EEE)_{\strict},$$
$$Q^{\op}\HOM(\ulG,\EEE)_{\strict} \isoto \HOM(\ulG,Q^{\op} \EEE)_{\strict}$$
natural in each variable. (Here we use the fact that the automorphisms in $Q^{\op}\EEE$ may be identified with the automorphisms in $\EEE$.)\\
(2) Taking the classifying space functor and composing the natural transformation$$ B\HOM(?,?) \Rightarrow \HHom_{\Top}(B(?),B(?)),$$
we get the maps 
$$BiQ^{\op}_{\bullet}\HOM(\ulG,\EEE)_{\strict} \to \HHom_{\Top}(BG,BiQ^{\op}_{\bullet}\EEE)$$
$$BQ^{op}\HOM(\ulG,\EEE)_{\strict} \to \HHom_{\Top}(BG,BQ^{\op}\EEE)$$
natural in each variable.\\
(3) Take the $n$-th homotopy group and the same argument in \textbf{Definition \ref{lax Sherman df}} and \textbf{Variant \ref{variant}}, for any connected pointed CW-complexes $X$, we get the following maps which shall be called the {\it{$n$-th original Sherman map}} (resp. {\it{$n$-th reduced original Sherman map}})
$$\Sh^{\ori}(X,\EEE)_n:R_n(X,\EEE) \to [B\pi_1(X),\Omega^{n+1}BQ^{\op}\EEE]$$
$$\tilde{\Sh}^{\ori}(X,\EEE)_n:\tilde{R}_n(X,\EEE) \to [B\pi_1(X),(\Omega^{n+1}BQ^{\op}\EEE)_0]_{\ast}$$
natural in each variables. Here we are using the $iS^{op}_{\bullet}=Q^{\op}$ theorem below. (c.f. \textbf{Review \ref{iS=Q}})
\end{ex}

\noindent
To compare the two Sherman maps, we will review the proof of $iS^{\op}_{\bullet}=Q^{\op}$ in \cite{Wal85}. (See also \cite{Gil81}.)

\begin{rev}[$iS^{\op}=Q^{\op}$ theorem] \label{iS=Q} $ $
\\
Let $\EEE$ be a small exact category.\\
(1) In \cite{Wal85} p.375, Waldhausen proved that the natural inclusion$\N Q^{\op} \EEE \hookrightarrow iQ_{\bullet}^{\op} \EEE$ induces a homotopy equivalence by the swallowing lemma (c.f. \cite{Wal85} p.352 Lemma 1.6.5.).\\ 
(2) We will define a natural morphism
$$iS^{\op,\edge}_{\bullet} \EEE \to iQ^{\op}_{\bullet} \EEE$$
as follows. For each non-negative integer $n$ and $X \in iS^{\op,\edge}_n \EEE$, we assign a sequence of $n$ composable morphisms in $Q^{\op}\EEE$,
$$\{ X(n+k/n-k+1) \rightarrowtail X(n+k/n-k) \twoheadleftarrow X(n+k+1/n-k) \}_{1 \leqq k \leqq n}.$$
This assignment defines the functor $iS^{\op,\edge}_n\EEE \to iQ^{\op}_n\EEE$ which is functorial in $n$ and full faithful and essentially surjective. Therefore we learn that it is a homotopy equivalence by \textbf{Lemma \ref{real lem}}.
\end{rev}

\begin{prop}[Compatibility] \label{compatible} $ $
\\
The Sherman map and the original Sherman map are compatible. Namely, if we consider an small exact category $\EEE$ as a category with fibrations and trivial equivalence relations, the following diagram is commutative for any pointed connected CW-complex $X$:
$$\xymatrix{
& [B\pi_1(X),\Omega^{n}\mathbb{K}^{\edge}(\EEE)] \ar[dd]^{\wr}\\
R_n(X,\EEE) \ar[ur]^{\!\!\!\Sh(X,\EEE)_n} \ar[dr]_{\Sh^{\ori}(X,\XXX)_n} \\
& [B\pi_1(X),\Omega^{n}BQ^{\op}\EEE]
}$$
where the vertical isomorphism is induced from the morphisms defined in \textbf{Review \ref{iS=Q}}. Similar statement is verified for reduced case.
\end{prop}

\begin{proof}
We shall only prove that the following diagrams are commutative for any group $G$:
$$\xymatrix{
BiS^{\op,\edge}_{\bullet}\HOM(\ulG,\XXX)_{\strict} \ar@{}[dr]|{\text{I}} \ar[r]^{\sim} \ar[d] &
B\HOM(\ulG,iS_{\bullet}^{\op,\edge}\XXX)_{\strict} \ar[d] \ar[r] \ar@{}[dr]|{\text{III}} &
\HHom_{\Top}(BG,BiS^{\op,\edge}_{\bullet}\XXX) \ar[d]\\
BiQ^{\op}_{\bullet}\HOM(\ulG,\XXX)_{\strict} \ar[r]^{\sim} \ar@{}[dr]|{\text{II}}&
B\HOM(\ulG,iQ_{\bullet}^{\op}\XXX)_{\strict} \ar[r] \ar@{}[dr]|{\text{IV}} &
\HHom_{\Top}(BG,BiQ^{\op}_{\bullet}\XXX)\\ 
BQ^{\op}\HOM(\ulG,\XXX)_{\strict} \ar[r]^{\sim} \ar[u] &
B\HOM(\ulG,Q^{\op}\XXX)_{\strict} \ar[u] \ar[r] &
\HHom_{\Top}(BG,BQ^{\op}\XXX) \ar[u]
}$$
where the vertical morphisms are defined in \textbf{Review \ref{iS=Q}}. The commutativity of I and II can be easily checked before taking the functor $B$. The commutativity of III and IV are consequence of the fact that 
$$B\HOM(?,?)_{\strict} \Rightarrow \HHom_{\Top}(B(?),B(?))$$
is a natural transformation by \textbf{Example \ref{Sherman argument}}.
\end{proof}

\noindent
Now we will state the universal property of algebraic $K$-theory associated with semi simple exact categories.

\begin{nt} $ $
\\
(1) We will let $\FFF$ denote the category of finite pointed connected CW-complexes.\\
(2) Consider the functors $F$ and $G$ from $\FFF$ to the category of commutative monoids, a natural transformation $t:F \to G$ is said to be {\it{universal}} if for any connected $H$-space $Z$ and any natural transformation $F \to [-,Z]_{\ast}$, the following diagram can be uniquely completed:
$$\xymatrix{
F \ar[r]^{t} \ar[dr] & G \ar@{-->}[d]\\
& [-,Z]_{\ast}
}$$
\end{nt}

\begin{thm}[\cite{She92} p.173 Corollary 5.2.] \label{universality} $ $
\\
If $\EEE$ is a semi-simple exact category, that is, all short exact sequences in $\EEE$ split, then the natural transformation 
$$\tilde{\Sh^{\ori}}(?,\EEE)_0:\tilde{R_0}(\pi_1(?),\EEE) \to [?,{(\Omega BQ^{\op}\EEE)}_0]_{\ast}$$
is universal. By \textbf{Proposition \ref{compatible}}, we learn that $\tilde{\Sh}(?,\EEE)_0$ is also universal. 
\end{thm}

\begin{his} $ $
\\
The original version of the statement above is first announced in \cite{Ger73}, and used to prove Gersten's conjecture for the case of discrete valuation rings whose residue fields are finite, and proved by Quillen in \cite{Hil81}. The variant statement is appeared in \cite{Lev97} and used to define the Adams operations on algebraic $K$-theory for (divisorial) schemes.
\end{his}

\section{Proof of the main theorem}

In this section, we will prove the main theorem by using the results in previous sections. The argument is divided the following a pair of lemmas.

\begin{lem} \label{key lemma 1} $ $
\\ 
To prove main theorem, we shall only check the following assertion:\\
For any group $G$, the canonical map induced from the inclusion map $\MMM(k) \hookrightarrow \AAA$,
$$R_0(G,\MMM(k)) \to R_0^{\lax}(G,\AAA)$$
is zero. Here $\AAA$ is defined in \textbf{Example \ref{ex of cat with er}} {\rm (5)}.
\end{lem}

\begin{proof}
First by the d\'evissage theorem \cite{Qui73} p.112 Theorem 4, the inclusion functor $\MMM(k) \hookrightarrow \MMM^1(R)$ induces isomorphisms on higher algebraic $K$-groups. So we may prove the inclusion functor $\MMM(k) \hookrightarrow \MMM(R)$ induces the zero maps on higher algebraic $K$-groups. We have the following commutative diagram for each $X \in \FFF$:
$$\xymatrix{
\tilde{R_0}(\pi_1(X),\MMM(k)) \ar[r] \ar[d]_{\tilde{\Sh}} & \tilde{R_0}(\pi_1(X),\MMM(R)) \ar[d]_{\tilde{\Sh}} & \tilde{R_0}(\pi_1(X),\PPP(R)) \ar[l] \ar[d]^{\tilde{\Sh}}\\
[X,{(\mathbb{K}(\MMM(k)))}_0]_{\ast} \ar[r] & [X,{(\mathbb{K}(\MMM(R)))}_0]_{\ast} & [X,{(\mathbb{K}(\PPP(R)))}_0]_{\ast} \ar[l]^{\sim}_{\text{I}}
}$$
where the morphism I is an isomorphism by the resolution theorem \cite{Qui73} p.108 Theorem 3.
It is well-known that $\mathbb{K}(\PPP(R))$ is a $H$-space (See for example \cite{Gra76}), $\MMM(k)$ is semi-simple and by the universal property \textbf{Theorem \ref{universality}}, we learn that we shall only prove the composition
$$\tilde{R_0}(\pi_1(X),\MMM(k)) \to \tilde{R_0}(\pi_1(X),\MMM(R)) \overset{\tilde{\Sh}}{\to} [X,{(\mathbb{K}(\MMM(R)))}_0]_{\ast}$$
is the zero map for any $X \in \FFF$. (The argument above is due to \cite{She82} p.240 and originally \cite{Ger73}.) Now remember that $\AAA$ is a category with fibrations and equivalence relations satisfying the cogluing axiom by \textbf{Example \ref{ex of cat with er}} (5), \textbf{Example \ref{cat with fib ex}} (3), and \textbf{Example \ref{ex of cogluing ax}} (2). Notice the inclusion functor factor through $\MMM(k) \hookrightarrow \AAA \hookrightarrow \MMM(R)$ and we have the following commutative diagram for each $X \in \FFF$:
$$\xymatrix{
& \tilde{R_0}(\pi_1(X),\MMM(R)) \ar[r]^{\tilde{\Sh}} & [X,{(\mathbb{K}^{\edge}(\MMM(R)))}_0]_{\ast}\\
\tilde{R_0}(\pi_1(X),\MMM(k)) \ar[r] & \tilde{R_0}(\pi_1(X),\AAA) \ar[r]^{\tilde{\Sh}} \ar[u] \ar[d] & [X,{(\mathbb{K}^{\edge}(\AAA))}_0]_{\ast} \ar[u] \ar[d]^{\text{II}}\\
& \tilde{R_0^{\lax}}(\pi_1(X),\AAA) \ar[r]_{\tilde{\Sh}^{\lax}} & [X,{(\mathbb{K}^{\lax,\edge}(\AAA))}_0]_{\ast}
}$$
where the morphism II is a injection by \textbf{Theorem \ref{rectification variant}}. Hence we get the result.
\end{proof}

\begin{lem} \label{key lemma 2} $ $
\\
The assertion in \textbf{Lemma \ref{key lemma 1}} is true.
\end{lem}

\begin{proof}
Let $G$ be a group and $(X,\rho_X)$ be a representation in $\HOM(\ulG,\MMM(k))_{\strict}$. Assume rank of $X$ as a $k$-vector space is $m$. Then there is a short exact sequence
\begin{equation} \label{crutial exact seq 2}
0 \to R^{\oplus m} \overset{\pi}{\to} R^{\oplus m} \overset{p}{\to} X \to 0.
\end{equation}
For each $g \in G$, we have a lifting of $\rho_X(g)$, that is, a $R$-module homomorphism $\tilde{\rho}(g):R^{\oplus m} \to R^{\oplus m}$ such that $\tilde{\rho}(g) \mod \pi=\rho_X(g)$. By Nakayama's lemma, we learn $\tilde{\rho}(g)$ is an isomorphism of $R$-modules. Obviously assignment $\tilde{\rho}:G \to \Aut(R^{\oplus m})$ defines a lax representation $(R^{\oplus m},\tilde{\rho})$ in $\LAX(\ulG,\AAA)_{\strict}$ and we have a short exact sequence
\begin{equation}  
(R^{\oplus m},\tilde{\rho}) \overset{\pi}{\to} (R^{\oplus m},\tilde{\rho}) \overset{p}{\to} (X,\rho_X)
\end{equation}
in $\LAX(\ulG,\AAA)_{\strict}$. Notice that proving $\pi$ is a strict deformation, we need the assumption that $R$ is commutative!! (Compare \textbf{Example \ref{noncom}}.) So we have an identity 
$$[(X,\rho_X)]=[(R^{\oplus m}, \tilde{\rho})]-[(R^{\oplus m},\tilde{\rho})]=0$$
in $R_0^{\lax}(G,\AAA)$. Hence we get the result. 
\end{proof}

\end{document}